\documentclass{amsart}
\usepackage{amscd,amssymb}
\usepackage{apn} 

\usepackage{graphicx}

\usepackage[margin=1.5in]{geometry}

\usepackage[pdftex]{hyperref}
\hypersetup{colorlinks,citecolor=blue,linktocpage}

\numberwithin{equation}{section}

\usepackage[all]{xy}

\usepackage[colorinlistoftodos]{todonotes}

\author{Valery Alexeev and Rita Pardini} 
\title{Non-normal abelian covers}

\begin{document}

\begin{abstract} 
  An abelian cover is a finite morphism $X\to Y$ of varieties which is
  the quotient map for a generically faithful action of a finite
  abelian group $G$. Abelian covers with $Y$ smooth and $X$ normal were
  studied in \cite{Pardini_AbelianCovers}.

  Here we study the non-normal case, assuming that $X$ and $Y$ are
  $S_2$ varieties that have at worst normal crossings
  outside a subset of codimension $\ge 2$. Special attention is paid
  to the case of $\Z_2^r$-covers of surfaces, which is used in
  \cite{AlexeevPardini_CB} to construct explicitly compactifications
  of some components of the moduli space of surfaces of general type.
\end{abstract}

\maketitle
\tableofcontents

\section*{Introduction}
An abelian cover is a finite morphism $X\to Y$ of varieties which is
the quotient map for a generically faithful action of a finite abelian
group $G$. This means that for every component $Y_i$ of $Y$ the
$G$-action on the restricted cover $X\times_Y Y_i\to Y_i$ is faithful.
The paper
\cite{Pardini_AbelianCovers} contains a comprehensive theory of such
covers in the case when $Y$ is smooth and $X$ is normal. The covers
are described in terms of the \emph{building data} consisting of
branch divisors $D_{H_i,\psi_i}$ ranging over cyclic subgroups
$H_i\subset G$, and line bundles $L_{\chi}$ with $\chi$ ranging over
the character group of $G$. This collection must satisfy the
\emph{fundamental relations}.

Here, we extend this theory to the case of singular varieties. Namely,
we allow $X$ and $Y$ to be varieties satisfying Serre's condition
$S_2$ and having double crossing singularities in codimension 1, which
we abbreviate to g.d.c. for ``generically double crossings''
(see \S\ref{ssec:Ysmooth} for the precise definition).
 Our
interest in this case lies in applications to the moduli theory. Such
non-normal abelian covers appear in our work \cite{AlexeevPardini_CB}
where we explicitly construct compactifications of moduli spaces of
some Campedelli and Burniat surfaces by adding stable surfaces on the
boundary. ``Stable surfaces'' here are in the sense of
\cite{KollarShepherdBarron}: they have slc (semi log canonical)
singularities and ample canonical class.

In this paper, we give a comprehensive treatment of the situation.  In
Section~\ref{ssec:Ysmooth} we show that the theory of standard covers
of \cite{Pardini_AbelianCovers} has a very natural extension to the
case when $Y$ is still smooth but $X$ is possibly g.d.c.. In
Section~\ref{ssec:Ynormal} we extend it to the case of normal base by
an $S_2$-fication trick.  In Section~\ref{ssec:Ynonnormal} we prove
that a cover with non normal $Y$ can be obtained by gluing a cover
over the normalization $\wY$, and we spell out which additional data
must be specified.

In Section~\ref{sec:geometry} we study the singularities of covers.
We determine the conditions for $X$ to have slc singularities,
to be Cohen-Macaulay, and we determine the index of the
canonical divisor in the situations appearing in common applications.

In Section~\ref{sec:Z2-covers} we treat in detail the special case
when the group $G$ is $\bZ_2^r$ and $\dim X=\dim Y=2$, as in
\cite{AlexeevPardini_CB}. We restrict ourselves to the situation where
the base $Y$ is smooth or has two smooth branches meeting transversally,
and the components of branch divisors and the double locus are smooth and
have distinct tangent directions at the points of intersection,
i.e. locally they look like a 
collection of lines in the plane. In this situation, we give a
complete classification of the covers and the singularities of
$X$. The answer is contained in nine tables. Some of these covers
appear on the boundary of moduli of Campedelli and Burniat surfaces,
but the full list is longer.

\begin{notations-nonum}
  $G$ denotes a finite abelian group. We work with equidimensional
  varieties defined over an algebraically closed field $\bK$ whose
  characteristic does not divide the order of $G$.  We denote by $G^*$
  the group $\Hom(G,\bK^*)$ of characters of $G$, and we write it
  multiplicatively. The abbreviations 
  \emph{lc} and {\em slc} stand for \emph{log canonical} and {\em semi
    log canonical}.  (cf.  \S \ref{sec:geometry} for the definitions).
  $\wX$, $\wC$, etc. denote the normalization of $X$, $C$, etc.  We
  use the additive and multiplicative notation for line bundles and
  divisors interchangeably. Linear equivalence will be denoted by
  $\sim$.
\end{notations-nonum}

\begin{acknowledgments} 
  The first author was partially supported by NSF under DMS 0901309.
  The second author wishes to thank Miles Reid and  Angelo
  Vistoli for several useful communications. We also thank the referee
  for many useful comments and corrections.

Part of this work was done while both authors were visiting MSRI in
the spring of 2009.  This project was partially supported by the Italian
PRIN 2008 project {\it Geometria delle variet\`a algebriche e dei loro
  spazi di moduli}.  The second author is a member of GNSAGA of INDAM.

\end{acknowledgments}

\section{General structure of abelian covers}

\subsection{Setup}\label{sec:setup}
We recall some basic facts about Serre's condition $S_2$ and the
$S_2$-fication of a coherent sheaf. For a comprehensive treatment, the
reader may consult \cite[5.9-11]{EGA4-2}, where the latter appears
under the name ``$Z^{(2)}$-closure''.

All varieties below are assumed to be reduced, equidimensional, but
possibly reducible. 
Let $\cF$ be a coherent sheaf on $X$ all of whose
associated components are irreducible components of $X$.  Then there
exists a unique \emph{$S_2$-fication}, or \emph{saturation in
  codimension 2}, a coherent sheaf defined by
\begin{displaymath}
  S_2(\cF)(V) = \varinjlim_{U\subset X,\ \codim (X\setminus U)\ge2}
  \cF(V\cap U)
\end{displaymath}

The sheaf $S_2(\cF)$ is $S_2$, and $\cF$ is $S_2$ iff the map $\cF\to
S_2(\cF)$ is an isomorphism. In particular, for $\cF=\cO_X$ one
obtains the $S_2$-fication $S_2(X) \to X$, which is dominated by the
normalization of~$X$.

On a normal variety $X$, an $S_2$-sheaf is the same as a \emph{reflexive
  sheaf}, satisfying $\cF^{**}=\cF$, see \cite{Bourbaki_AC7}. Further,
reflexive sheaves of rank 1 are the same as \emph{divisorial sheaves},
isomorphic to $\cO_X(D)$ for some Weil divisor $D$, see
e.g. \cite[App.to \S1]{Reid_Can3folds}.
On a smooth (or factorial) variety Weil divisors are the same as
Cartier divisors, and rank 1 $S_2$ sheaves are the same as invertible
sheaves. 

Let $G$ be a finite abelian group. An {\em abelian cover} with Galois
group $G$, or {\em $G$-cover}, is a finite morphism $X\to Y$ of
varieties which is the quotient map for a generically faithful action
of a finite abelian group $G$. This means that for every component
$Y_i$ of $Y$ the $G$-action on the restricted cover $X\times_Y Y_i\to
Y_i$ is faithful.  An {\em isomorphism} of $G$-covers $\pi_1\colon
X_1\to Y$, $\pi_2\colon X_2\to Y$ is an isomorphism $\phi\colon X_1\to
X_2$ such that $\pi_1=\pi_2\circ \phi$.

The $G$-action on $X$ with $X/G=Y$ is equivalent to a decomposition:
\begin{equation}\label{eq:decomposition}
  \pi_*\OO_X=\bigoplus_{\chi\in G^*}\cF_{\chi},\qquad \cF_{1}=\cO_Y
\end{equation}
where $G$ acts on $\cF_{\chi}$ via the character $\chi$.  If $\pi$ is
Galois then each $\cF_{\chi}$ has rank 1: if $y\in Y$ is a general
closed point, then $G$ acts freely on $\pi\inv (y)$, so it acts on
$\OO_{\pi\inv(y)}=\oplus_{\chi}(\cF_{\chi}\otimes \K(y)) $ as the
regular representation. Thus, $\cF_{\chi}\otimes \K(y)$ is
$1$-dimensional for every $\chi$. When the sheaves $\cF_{\chi}$ are
locally free, it is customary to write $\cF_{\chi}=L_{\chi}\inv$, with
$L_{\chi}$ a line bundle.

\begin{lemma}\label{lem:flatness-and-CM}
  \begin{enumerate}
  \item The sheaf $\cO_X$ is $S_n$ for some $n$ iff every
    $\cF_{\chi}$ is $S_n$. 
  \item If $\pi\colon X\to Y$ is flat then $X$ is CM iff $Y$ is CM.
  \item If $Y$ is smooth and $X$ is $S_2$ then $\pi$ is flat and $X$
    is CM.
  \end{enumerate}
\end{lemma}
\begin{proof}
  (1) is clear by definition of depth. 

  (2) $\pi$ is flat iff every $\cO_Y$-module $\cF_{\chi}$ is
  invertible. Then each $\cF_{\chi}$ is CM iff $\cO_Y$ is. 

  (3) On a smooth variety every divisorial sheaf is invertible, and so
  flat. Now (2) applies. 
\end{proof}

A $G$-cover $\pi\colon X\to Y$, where $X$ and $Y$ are $S_2$ varieties,
is determined by its restriction to the complement of a closed subset
of codimension $\ge 2$:
\begin{lemma}\label{lem:remove-codim2}
  Let $Y$ be an $S_2$ variety, $Y_0\subseteq Y$ an open subset with
  $\codim (Y\setminus Y_0)\ge2$, and $\pi_0\colon X_0\to Y_0$ a
  $G$-cover with $X_0$ an $S_2$ variety. Then there exist a unique
  $S_2$ variety $X$ and a $G$-cover $\pi\colon X\to Y$ whose
  restriction to $Y_0$ is $\pi_0$.
\end{lemma}
\begin{proof}
  For the existence, we take $\cO_X:= i_* \cO_{X_0}$, where $i\colon
  Y_0\to Y$ is the inclusion. Then $\cO_X=\oplus_{\chi\in G^*} \cF_{\chi}$, where each
  $\cF_{\chi}$ is a rank 1 $S_2$-sheaf. The algebra structure on $\cO_X$
  is defined as follows. For an open set $U\subset X$ and sections
  $s\in \cF_{\chi}(U)$, $s'\in \cF_{\chi'}(U)$, their product is
  \begin{displaymath}
    s|_{U\cap X_0} \cdot s'|_{U\cap X_0} \in \cF_{\chi\chi'}(U\cap X_0)
    = \cF_{\chi\chi'}(U),
  \end{displaymath}
  since $\codim_U(U\setminus U\cap X_0)\ge2$ and $\cF_{\chi}$ is saturated in
  codimension 2. Thus, $X:=\Spec_{\cO_Y}\cO_X$ is an $S_2$ variety with a
  finite morphism to $Y$. 
  The $G^*$-grading on $\cO_X$ defines the $G$-action on $X$.
  By construction, the eigenspace $\cF_1$ for the trivial character is
  $i_*\cO_{Y_0} = \cO_Y$. Therefore, $X/G=Y$.

  Uniqueness follows from the uniqueness of the $S_2$-fication.
\end{proof}
Given a $G$-cover $\pi\colon X\to Y$ and
 an irreducible  subset $S\subset Y$, we
define the {\em inertia subgroup} $H_S$ of $S$ to be the subgroup of
$G$ consisting of the elements that fix $\pi\inv(S)$ pointwise, or,
equivalently since $G$ is abelian, 
that fix an irreducible component of $\pi\inv(S)$
pointwise.  
The {\em branch locus} $D_{\pi}$ of $\pi$ is the set of
points of $Y$ whose inertia subgroup is not trivial (notice that we
regard $D_{\pi}$ simply as a set, without giving it a scheme
structure).  If $\pi$ is flat, then $D_{\pi}$ is a
divisor by \cite[Thm. 6.8]{AltmanKleiman_GrothDuality}. 
If $F$ is an irreducible divisor of
$Y$ such that $X$ is generically smooth along $\pi\inv(F)$, then the
natural representation $\psi$ of $H_F$ on the tangent space $T_{X,R}$
at the generic point of an irreducible component $R$ of $\pi\inv(F)$
is faithful, hence $H_{F}$ is cyclic (cf.  \cite[\S
1]{Pardini_AbelianCovers}). Notice that $\psi$ does not depend on the
choice of the component $R$ of $\pi\inv(F)$ since $G$ is abelian.

\subsection{Standard covers}\label{ssec:standard}

In this section we recall, in a form which is convenient for our later
applications, the definition of standard abelian covers, a class of
flat abelian covers that can be constructed from a collection of line
bundles and effective divisors on the target variety
(cf. \cite{Pardini_AbelianCovers}, \cite{FantechiPardini}).  The
prototypical example is the classical construction of a double cover
of a variety $Y$ from the data of an effective divisor $D$ on $Y$ and
a line bundle $L$ such that $2L\sim D$.  \medskip

Let $Y$ be a variety. A set of {\em building data for a standard
  $G$-cover} $\pi\colon X\to Y$ consists of the following:
\begin{itemize}
\item irreducible effective Cartier divisors $D_1,\dots D_k$ (possibly not distinct),
\item for each $D_i$ a pair $(H_i, \psi_i)$, where $H_i$ is a cyclic
  subgroup of $G$ of order $m_i$ and $\psi_i$ is a generator of the
  group of characters $H_i^*$,
\item line bundles $L_{\chi}$, for $\chi\in G^*\setminus\{1\}$.
\end{itemize}
Moreover we assume that these data satisfy the so called {\em
  fundamental relations}:
\begin{equation}\label{eq:fundrel}
  \forall \chi,\chi', \quad
  L_{\chi}+ L_{\chi'} \sim  L_{\chi\chi'} + 
  \sum_i\epsi^i_{\chi,\chi'} D_i,
\end{equation}
where for a character $\chi$ we write
$\chi|_{H_i}=\psi_i^{a_{\chi}^i}$, with $0\le a^i_{\chi}< m_i$, and we
define
$\epsi^i_{\chi,\chi'}:=[\frac{a_{\chi}^i+a_{\chi'}^i}{m_i}]$. Observe
that $\epsi^i_{\chi, \chi'} $ is equal either to $0$ or to $1$.

We call the divisors $D_i$, together with the pairs $(H_i,\psi_i)$,
the {\em branch data} of the cover.  An equivalent way of describing
the branch data, and therefore the building data, is to give for each
pair $(H,\psi)$, with $H\subset G$ a cyclic subgroup and $\psi\in H^*$
a generator, the divisor
$D_{H,\psi}=\sum_{\{i|(H_i,\psi_i)=(H,\psi)\}}D_i$. This is the
notation used in \cite{Pardini_AbelianCovers}.

\begin{remark}\label{rem:branchdata}
  If the group $\Pic(Y)$ has no $m$-torsion, where $m=|G|$, then the
  branch data determine the building data by
  \cite[Prop. 2.1]{Pardini_AbelianCovers}. In general, the branch data
  are enough to determine the local geometry of the cover
  (cf. Proposition \ref{prop:bdata}, (2)).
\end{remark}
\begin{remark}
  When $G=\Z_2^r$, it is enough to associate with every divisor $D_i$
  a nonzero element $g_i\in G$, the generator of $H_i$. Also, the
  definition of $\epsi^i_{\chi,\chi'} $ is simpler: $\epsi^i_{\chi,
    \chi'} $ is equal to 1 if $\chi(g_i)=\chi'(g_i)=-1$ and it is
  equal to 0 otherwise.
\end{remark}
We now explain how to construct a $G$-cover from a set of building
data.  Choose $\chi_1,\dots \chi_s\in G^*$ such that $G^*$ is the
direct sum of the cyclic subgroups generated by the $\chi_j$. Denote
by $d_j$ the order of $\chi_j$ and write $L_j:=L_{\chi_j}$ and
$a^i_j:=a^i_{\chi_j}$. By \cite[Prop.2.1]{Pardini_AbelianCovers} for
$j=1,\dots s$ there exist isomorphisms:
$$\fie_j: L_j^{\otimes d_j}\iso \OO_Y(\sum_i \frac{d_ja^i_j}{m_i}D_i).$$
Notice that the coefficients $\frac{d_ja^i_j}{m_i}$ in the above formula are integers. 
Using formulae (2.15) of \cite{Pardini_AbelianCovers} and the
isomorphisms $\fie_j$ above, 
one constructs for each pair $\chi,\chi'$
of non trivial characters an isomorphism $$\fie_{\chi,\chi'}\colon
L_{\chi}\inv \otimes L_{\chi'}\inv \iso L_{\chi\chi'}\inv (-\sum
\epsi^i_{\chi,\chi'}D_i)$$ such that for every $\chi, \chi',\chi''\in
G^*$ the following diagram commutes (we set $L_1=\OO_Y$):
\begin{equation}\label{diag:mu_chi}
  \begin{CD}
    L_{\chi}\inv \otimes L_{\chi'}\inv \otimes L_{\chi''}\inv  @>>>L_{\chi\chi'}\inv (-\sum_i\epsi^i_{\chi,\chi'}D_i)\otimes L_{\chi''}\inv  \\
    @VVV @VVV\\
    L_{\chi}\inv \otimes L_{\chi'\chi''}\inv
    (-\sum_i\epsi^i_{\chi',\chi''}D_i)@>>>L_{\chi\chi'\chi''}\inv
    (-\sum_i\de^i_{\chi,\chi',\chi''}D_i)
  \end{CD}
\end{equation}
where
$\de^i_{\chi,\chi',\chi''}=\epsi^i_{\chi\chi',\chi''}+\epsi^i_{\chi,\chi'}=\epsi^i_{\chi,
  \chi'\chi''}+\epsi^i_{\chi',\chi''}$ and the maps are induced by the
$\fie_{\chi,\chi'}$ in the obvious way.  We denote by
$\mu_{\chi,\chi'}\colon L_{\chi}\inv\otimes L_{\chi'}\inv\to
L_{\chi\chi'}\inv$ the maps induced by composing $\fie_{\chi,\chi'}$
with the inclusion $L_{\chi\chi'}\inv (-\sum
\epsi^i_{\chi,\chi'}D_i)\hookrightarrow L_{\chi\chi'}\inv$. By the
commutativity of diagram \eqref{diag:mu_chi}, the collection of maps
$\mu_{\chi,\chi'}$ defines on $\cE:=\OO_Y\oplus \bigoplus_{\chi\ne
  1}L_{\chi}\inv$ a commutative and associative algebra structure
compatible with the $G$-action defined by letting $G$ act trivially on
$L_1=\cO_Y$ and via the character $\chi$ on $L_{\chi}\inv$ for
$\chi\ne 1$. We define $X:=\Spec \cE$ with the natural map $\pi \colon
X\to Y$ to be a {\em standard $G$-cover} associated with the given set
of building data.  Notice that, since the $L_{\chi}\inv$ are
locally free, $\pi$ is flat and $X$ is $S_2$ if $Y$ is.

$X$ can be described locally above a point $y\in Y$ as follows. Up to
shrinking $Y$, we may assume that all the $L_{\chi}$ are trivial and
that the $D_i$ are defined by equations $\si_i$. If we denote by
$z_{\chi}$ a coordinate on $L_{\chi}\inv$, $\chi\in
G^*\setminus\{1\}$, then $X$ is given inside the vector bundle
$V(\oplus_{\chi\ne 1}L_{\chi}\inv)\cong Y\times \bK^{m-1}$ by the
following set of equations:
\begin{equation}\label{eq:local}
  z_{\chi}z_{\chi'}=c_{\chi,\chi'}\Pi_1^k \si_i^{\epsi^i_{\chi,\chi'}} z_{\chi\chi'}, \quad \chi,\chi'\in G^*\setminus\{1\},
\end{equation}
where the $c_{\chi,\chi'}$ are nowhere vanishing regular functions and
for $\chi=1$ we set $z_{\chi}=1$.  For $1\ne \chi\in G^*$, denote by
$d$ the order of $\chi$ and write $\chi|_{H_i}=\psi_i^{a_i}$, with
$0\le a_i<m_i:=|H_i|$. Eliminating between the equations
\eqref{eq:local}, one gets 
\begin{equation}\label{eq:powerz}
  z_{\chi}^d=b_{\chi}\Pi_1^k\si_i^{\frac{da_i}{m_i}},
\end{equation}
where $b_{\chi}$ is a nowhere vanishing function.  It follows
immediately that $X$ is a variety: indeed, using the decomposition of
$\pi_*\OO_X$ into $G$-eigenspaces, we may assume  that a nilpotent
element is locally of the form $fz_{\chi}$ for some character $\chi$
and some regular function $f$. 
Then by
\eqref{eq:powerz},
 $(fz_{\chi})^k=0$ for some $k$ only if $f=0$.
Using the local equations \eqref{eq:local}, one can also show the
following:
\begin{lem}\label{lem:inertia} Notation as above. Let $\pi\colon X\to
  Y$ be a standard $G$-cover and $y\in Y$ be a point. The inertia
  subgroup $H_y$ of $y$ is equal to $\sum_{\{i|y\in D_i\}}H_i$.
\end{lem}
\begin{proof} Since the question is local on $Y$, we may assume that
  $X$ is given by the equations \eqref{eq:local}.  Let $x\in X$ be a
  point lying above $y$. Then by \eqref{eq:powerz} the coordinate
  $z_{\chi}(x)$ does not vanish iff $\chi|_{H_i}=1$ for every $i$ such
  that $y\in D_i$.  Since an element $g\in G$ fixes $x$ if and only if
  for every $\chi\in G^*$ such that $\chi(g)\ne 1$ the coordinate
  $z_{\chi}(x)$ vanishes, this remark proves the claim.
\end{proof}
 
Given a set of building data, the construction of the standard
$G$-cover $\pi\colon X\to Y$ depends of course on the choice of the
characters $\chi_1, \dots \chi_s$ and of the isomorphims
$\fie_j$. Assume that $\chi'_1,\dots \chi'_t$ are another set of
characters of $G$ such that $G^*$ is the direct sum of the cyclic
subgroups generated by the $\chi'_l$. Let $d'_l$ be the order of
$\chi'_l$, $i=1,\dots t$; then by \eqref{eq:powerz} the multiplication
maps \underline{induce} for $l=1,\dots t$ isomorphisms $\fie'_l\colon
L_{\chi'_l}^{\otimes d'_l}\iso \OO_Y(\sum_i \frac{k_lb^i_l}{m_i}D_i)$,
where $0\le b^i_l<m_i$ and $\chi'_l|_{H_i}=\psi_i^{b^i_l}$.  By the
associativity and commutativity of the multiplication the algebra
structure defined on $\OO_Y\oplus \bigoplus_{\chi\ne 1}L_{\chi}\inv$
by the $\fie'_l$ is the same as that induced by the $\fie_j$.  Hence
it is enough to analyze to what extent the isomorphism class of $\pi$
depends on the~$\fie_j$:
\begin{proposition}\label{prop:bdata}
  \begin{enumerate}
  \item {\em (Global case).} If $H^0(\OO_Y^*)=\bK^*$, then the
    building data determine $\pi\colon X\to Y$ up to isomorphism of
    $G$-covers.
  \item In general, given two standard covers $\pi_i\colon X_i\to Y$,
    $i=1,2$, with the same building data, there exists an \'etale
    cover $Y'\to Y$ such that, after base change with $Y'\to Y$,
    $\pi_1$ and $\pi_2$ give isomorphic $G$-covers.
  \end{enumerate}
\end{proposition}
\begin{proof} (2) We use the notation introduced above.  Let $\cE$, $\cE'$
  be two $\OO_Y$-algebra structures on $\OO_Y\oplus \bigoplus_{\chi\ne
    1}L_{\chi}\inv$ given by isomorphisms $\fie_j$, respectively
  $\fie'_j$.  The isomorphisms $\fie_j$, $\fie'_j$ differ by an
  automorphism of $L_j^{\otimes d_j}$, namely by multiplication by an
  element $k_j\in H^0(\OO_Y^*)$. This automorphism is induced by an
  automorphism of $L_j$ iff $k_j$ has a $d_j$-th root $h_j\in
  H^0(\OO_Y^*)$. So, up to taking an \'etale cover, one can assume
  that the roots $h_j$ exist. By formulae (2.15) of
  \cite{Pardini_AbelianCovers}, the $h_j$ can be used to define for
  all $\chi\in G^*\setminus\{1\}$ automorphisms $\psi_{\chi}$ of
  $L_{\chi}\inv $ that commute with the isomorphisms
  $\fie_{\chi,\chi'}$ and $\fie'_{\chi,\chi'}$.

  To prove statement (1), just observe that if $H^0(\OO_Y^*)=\bK^*$ no
  base change is necessary to construct the isomorphism above.
\end{proof}
\begin{remark}\label{rem:local} Let $\pi\colon X\to Y$ be a $G$-cover
  with branch data $D_i, (G_i,\psi_i)$, let $y\in Y$ and let $\si_i$
  be local equations for $D_i$ near $y$.  Combining Proposition
  \ref{prop:bdata} with the local equations \eqref{eq:local}, we see
  that, up to passing to an \'etale cover of $(Y,y)$, $X$ is defined
  locally near $y$ by the equations:
  \begin{equation}\label{eq:localc1}
    z_{\chi}z_{\chi'}=\prod_{i=1}^k \si_i^{\epsi^i_{\chi,\chi'}} z_{\chi\chi'}, \quad \chi,\chi'\in G^*\setminus\{1\},
  \end{equation}
\end{remark}

\subsection{Covers of smooth varieties}\label{ssec:Ysmooth}
Here we find conditions for a $G$-cover of a smooth variety to be
standard.  We keep the notation of the previous section.

\begin{definition}\label{defn:hurwitz}
  Let $Y$ be a smooth variety and let $\pi\colon X\to Y$ be a standard
  $G$-cover with building data $L_{\chi}$, $D_i$, $(H_i,\psi_i)$. By
  Lemma \ref{lem:inertia} the branch locus $D_{\pi}$ of $\pi$ is the
  support of the divisor $\sum_iD_i$.

  We define the {\em Hurwitz divisor} of $\pi$ as the $\Q$-divisor
  $D:=\sum_i\frac{m_i-1}{m_i}D_i$.  Notice that the support of $D$ is
  equal to $D_{\pi}$.
\end{definition}
 
We say that a variety is {\em d.c.} (has {\em double crossings}) if
every point is either smooth or analytically isomorphic to $xy=0$. We
say that a variety is
{\em g.d.c.} (has {\em generically double
  crossing}s) if it is d.c. outside a closed subset of codimension
$\ge 2$.

The following result
  generalizes the main result of \cite{Pardini_AbelianCovers}:

\begin{theorem}\label{thm:structure} Let $\pi\colon X\to Y$ be a
  $G$-cover such that $Y$ is smooth and $X$ is $S_2$. Then:
  \begin{enumerate}
  \item $X$ is normal iff $\pi$ is standard and every component of the
    Hurwitz divisor $D$ has multiplicity $<1$. 
  \item Assume that $\pi$ is standard.  Then $X$ is g.d.c. iff every
    component of $D$ has multiplicity $\le 1$.
  \item Assume that $X$ is g.d.c.. Then $\pi$ is standard iff for
    every irreducible divisor $F$ of $Y$ such that $X$ is singular
    above $F$ one has $H_{F}=\Z_2^s$ for some $s$.
  \end{enumerate}
\end{theorem}
\medskip

In the case $G=\Z_2^r$, which is of special interest to us because of
the applications in \cite{AlexeevPardini_CB}, Theorem \ref{thm:structure}
reads:
\begin{cor}\label{cor:structureZ2} Let $\pi\colon X\to Y$ be a
  $\Z_2^r$-cover such that $Y$ is smooth and $X$ is $S_2$. Then:
  \begin{enumerate}
  \item $X$ is normal iff $\pi$ is standard and every component of $D$
    has multiplicity $<1$.
  \item $X$ is g.d.c. iff $\pi$ is standard and every component of $D$
    has multiplicity $\le 1$.
  \end{enumerate}
\end{cor}

\begin{remark} Let $\pi\colon X\to Y$ be a standard $G$-cover with $Y$
  smooth and $X$ g.d.c.  and let $F$ be a component of the branch
  divisor $D_{\pi}$. By Lemma \ref{lem:inertia}, we have
  $H_{F}=\sum_{\{i|D_i=F\}}H_i$.  The pairs (subgroup, character)
  corresponding to $F$ can be determined as follows:
  \begin{itemize}
  \item Assume that $F$ has multiplicity $<1$ in the Hurwitz divisor
    $D$.  Then there is precisely one index $i$ with $D_i=F$.  In this
    case, $H_i=H_{F}$ and the character $\psi_i$ is given by the
    action of $H_i$ on the tangent space to $X$ at the generic point
    of an irreducible component of $\pi\inv(F)$
    (cf. \cite{Pardini_AbelianCovers}, \S 1 and \S 2).
  \item Assume that $F$ has multiplicity $=1$ in $D$.  Then there are
    precisely two indices $i_1$ and $i_2$ such that $D_{i_1}=D_{i_2}=
    F$ and $H_{i_1}$ and $H_{i_2}$ have order 2. So either
    $H_{F}=H_{i_1}=H_{i_2}$ or $H_{F}=H_{i_1}\oplus H_{i_2}$. In the
    latter case the proof of Theorem \ref{thm:structure} shows that
    $H_{i_1}$ and $H_{i_2}$ are generated by the elements of $H_{F}$
    that interchange the two branches of $X$ at a general point of
    $\pi\inv(F)$.
  \end{itemize}
\end{remark}
\begin{proof}[Proof of Theorem \ref{thm:structure}] Statement (1) is
  \cite{Pardini_AbelianCovers}, Thm. 2.1 and Cor.3.1.

  So consider the non-normal case. The cover $\pi$ is flat since $Y$
  is smooth and $X$ is $S_2$, hence we write as usual
  $\pi_*\OO_X=\OO_Y\oplus\bigoplus_{\chi\ne 1}L_{\chi}\inv$. The cover
  is standard if and only if there exist branch data $D_i$,
  $(H_i,\psi_i)$ such that for every $\chi,\chi'\in G^*\setminus\{1\}$
  the zero divisor of the multiplication map $\mu_{\chi,\chi'}\colon
  L_{\chi}\inv\otimes L_{\chi'}\inv\to L_{\chi\chi'}\inv$ is equal to
  $\sum_i\epsi^i_{\chi,\chi'}D_i$, where the $\epsi^i_{\chi,\chi'}$
  are defined in \S \ref{ssec:standard}.

  Notice that $X$, being $S_2$, is non-normal if and only if it is
  singular in codimension 1.  Fix a component $F$ of $D$ such that $X$
  is singular above $F$. Write $H:=H_{F}$. The cover $\pi$ factors as
  $X \to X/H\to Y$ and $F$ is not contained in the branch locus of the
  map $X/H\to Y$, hence $X/H$ is generically smooth over $F$. It
  follows that there is an element of $H$ that exchanges the two
  branches of $X$ at a general point of $\pi\inv(F)$.

  Let $\wX\to X$ be the normalization, let $\pi\unu\colon \wX \to
  Y$ be the induced $G$-cover, 
  let $(H', \psi')$ be the pair
  (subgroup, character) corresponding to $F$ for the cover $\pi\unu$
  and let $m'$ be the order of $H'$ (if $\pi\unu$ is not branched on
  $F$, we take $H'$ and $\psi'$ to be trivial).  Since the
  normalization map $\wX\to X$ is $G$-equivariant, we have a short
  exact sequence:
  \begin{equation}
    0\to H'\to H \to \Z_2\to 0.
  \end{equation}

  We consider the $H$-covers $p\colon X\to Z:=X/H$ and $p\unu\colon
  \wX\to \wX/H=Z$ and we study the algebras $\mathcal
  A:=p_*\OO_{X,F'}$ and $\cA\unu:=p\unu_*\OO_{\wX,F'}$, where $F'$
  is an irreducible component of the inverse image of $F$ in $Z$. We
  denote by $t\in \OO_{Z,F'}$ a local parameter.

  We distinguish three cases:

  \smallskip
  \noindent \underline{Case (a): $|H|=2$.}\par In this case
  $H'=\{0\}$, and $X$ is given locally by $z^2=at^2$, where $a\in
  \OO_{Z,F'}^*$ .  \smallskip

  \smallskip
  \noindent \underline{Case (b): $H$ is cyclic of order $2m'\ge
    4$.} \par
  Let $\psi\in H^*$ be a generator that restricts to $\psi'$ on
  $H'$. The algebra $\cA\unu$ is generated by elements $z, w$ such
  that:
  \begin{equation}
    z^{m'}=atw, \quad w^2=b
  \end{equation}
  where $a, b\in \OO_{Z, F'}^*$ and $H$ acts on $z$ via the character
  $\psi$ and on $w$ via the character $\psi^{m'}$.  The eigenspace
  corresponding to $\psi^j$ is generated by $z_j:=z^j$ for $0\le
  j<m'$, and by $z_j:=wz^{j-m'}$ for $m'\le j< 2m'$. Since the
  inclusion $\cA\subset \cA\unu$ is $G$-equivariant, ${\mathcal A}$ is
  generated by elements of the form $t^{a_j}z_j$ for suitable $a_j\ge
  0$.

  Since $H$ fixes $p\inv (F')$ pointwise, by the argument in the proof
  of Lemma \ref{lem:inertia} $\mathcal A$ is contained in the
  subalgebra $\mathcal B$ of $\cA\unu$ generated by
$$1, z^{m'}=tw, z_j,\quad 1\le j\le 2{m'}-1, \quad j\ne m'.$$ 
${\mathcal B}$ is also generated by $z_1=z, z_{m'+1}=wz$, with the
only relation $bz_1^2=z_{m'+1}^2$, hence $\Spec {\mathcal B}$ is
g.d.c. and the map $\Spec \cB\to \Spec\cA$ is an isomorphism.  So
${\mathcal A}=\mathcal B$.  

\smallskip
\noindent \underline{Case (c): $H$ is not cyclic.} \par
In this case $m'$ is even and $H\cong H'\times \Z_2$.  We denote by
$\psi\in H^*$ a character that restricts to $\psi'$ on $H'$ and by
$\phi$ the character such that $H'=\ker \phi$.  $\cA\unu$ is generated
by $z, w$ such that:
\begin{equation}
  z^{m'}=at, w^2=b,
\end{equation}
where $a, b\in \OO_{Z, F'}^*$ and $H$ acts on $z$ via the character
$\psi$ and on $w$ via the character $\phi$.  Arguing as in the
previous case, one checks that ${\mathcal A}$ is generated by:
$$1, z_1:=z,\dots, z^{m'-1}, tw,  z_{m'+1}:=zw, \dots, z^{m'-1}w.$$
$\mathcal A$ can also be generated by $z_1, z_{m'+1}$ with the only
relation $bz_1^2=z_{m'+1}^2$.  \smallskip

For $\chi_1,\chi_2\in G^*\setminus\{1\}$, denote by
$\epsi_{\chi_1,\chi_2}$ the order of vanishing on $F$ of the
multiplication map $\mu_{\chi_1,\chi_2}\colon L_{\chi_1}\inv \otimes
L_{\chi_2}\inv \to L_{\chi_1\chi_2}\inv$. Using the above analysis and
arguing as in the proof of \cite[Thm. 2.1]{Pardini_AbelianCovers}, one
obtains the following rules, up to exchanging $\chi_1$ and $\chi_2$:
\begin{itemize}

\item[Case (a):] $\epsi_{\chi_1,\chi_2}=2$ if $\chi_1,\chi_2\notin
  H^{\perp}$, \par\noindent $\epsi_{\chi_1,\chi_2}=0$ otherwise.

\item[Case (b):]

  For $i=1,2$, write $\chi_i|_H=\psi^{\al_im'+\be_i}$, where $\al_i=0$
  or $1$ and $0\le \be_i<m'$. Then:

  $\epsi_{\chi_1,\chi_2}=2$ if $\al_1=\al_2=1$, $\be_1=\be_2=0$,

  $\epsi_{\chi_1,\chi_2}=1$ if $\al_1=1$, $\be_1=0$, $\be_2>0$,

  $\epsi_{\chi_1,\chi_2}= [(\be_1+\be_2-1)/m']$ in the remaining
  cases.

\item[Case (c):] For $i=1,2$, write
  $\chi_i|_H=\phi^{\al_i}\psi^{\be_i}$,where $\al_i=0$ or $1$ and
  $0\le \be_i<m'$. Then:

  $\epsi_{\chi_1,\chi_2}=2$ if $\al_1=\al_2=1$, $\be_1=\be_2=0$,

  $\epsi_{\chi_1,\chi_2}=1$ if $\al_1=1$, $\be_1=0$, $\be_2>0$

  $\epsi_{\chi_1,\chi_2}= [(\be_1+\be_2)/m']$ in the remaining cases.

\end{itemize}

In the above analysis the group $\bZ_2^s$ appears in case (a) and case
(c) for $m'=2$.  In case (a), the cover $\pi$ is standard: $F$ appears
twice among the branch data, both times with label $H$.  In case (c),
$\pi$ is standard for $m'=2$: $F$ appears twice among the branch data,
with labels ${H_1}$ and ${H_2}$ corresponding to the subgroups of
order 2 of $H$ distinct from $H'$.  Moreover, it is not difficult to
check that in case (b) and in case (c) for $m'\ne 2$ the cover is not
standard. So we have proven (3) and also that every component of the
Hurwitz divisor $D$ of a standard g.d.c. cover has multiplicity $\le
1$.
  
Vice versa, assume that $\pi$ is standard and $F$ appears in $D$ with
multiplicity $\le 1$. If the multiplicity is $<1$ then the cover is
normal over $F$. If the multiplicity is equal to $1$, then $F$ appears
twice among the branch data, and the corresponding subgroups $H_1$ and
$H_2$ have order 2.  If $H_1=H_2$, then the cover is given over the
generic point of $F$ by the equation $z^2=ut^2$, with $u$ a unit; so
it is g.d.c.  If $H_1\ne H_2$, then the cover is given by the
equations $z_1^2=at$, $z_2^2=bt$, with $a$ and $b$ units. These equations  are equivalent to
$az_2^2=bz_1^2$, so the cover is  g.d.c.. This completes the proof of (2).
\end{proof}

\subsection{Covers of normal varieties}
\label{ssec:Ynormal}

Let $\pi\colon X\to Y$ be a $G$-cover such that $Y$ is normal and $X$ is  $S_2$. Let $Y_0$ be the
nonsingular locus of $Y$. Then the restriction $\pi_0\colon X_0\to Y_0$
is a $G$-cover, and by Lemma~\ref{lem:remove-codim2} $\pi$ is the unique
$S_2$-extension of $\pi_0$ to $Y$. Thus the theory in the normal case is
the immediate extension of the nonsingular case. We record the changes:
\begin{enumerate}
\item The sheaves $\cF_{\chi}$ are no longer invertible but they are
  $S_2$, i.e. in this case reflexive, divisorial sheaves. The
  multiplication maps are 
  \begin{displaymath}
    \cF_{\chi}\times \cF_{\chi'} \to
    \cF_{\chi}\otimes \cF_{\chi'} \to
    (\cF_{\chi}\otimes \cF_{\chi'})^{**} \to
    \cF_{\chi\chi'}.  \end{displaymath}

\item The branch divisors $D_g$ are Weil divisors.
\end{enumerate} 
Otherwise, the same fundamental relations between $\cF_{\chi}$ and
$D_g$ must hold. 

One has to be careful that the morphism $\pi$ may be not flat; indeed,
it is flat iff all $\cF_{\chi}$ are invertible. Also, for a singular
$Y$ the branch locus may have non-divisorial components.

\begin{example}
  Let $X=\bA^2=\Spec k[x,y]$, $G=\bZ_2$ acting by $x\mapsto -x$,
  $y\mapsto -y$, and let $Y$ be the quotient $\Spec k[x^2,xy,y^2]$, a
  quadratic cone. Then $\pi$ is ramified only over the vertex $P$ of the
  cone. The divisors $D_g$ are zero. The eigensheaves are
  $\cF_1=\cO_Y$ and $\cF_{-1}$, the divisorial sheaf corresponding to
  a line $\ell$ through the vertex. $\cF_{-1}$ is also isomorphic to
  the $\cO_Y$-submodule of $\cO_X$ generated by $x$ and $y$. 

  The fundamental relation in this case is $2\cF_{-1}=0$.
\end{example}

\subsection{Covers of non normal varieties}
\label{ssec:Ynonnormal}

Now we assume that $Y$ is a non normal g.d.c. and $S_2$ variety.  Let
$C$ be the divisorial part of the singular locus of $Y$, let
$\nu\colon\wY\to Y$ be the normalization, let $C'$ be the inverse
image of $C$ in $\wY$ and let $\wt{C'}\to C'$ be the normalization.
Since $Y$ is g.d.c., there is a biregular involution $\iota$ on
$\wt{C'}$ induced by the degree 2 map $\wt{C'}\to C'\to C$.  (If the
components of $Y$ are smooth, then $\wt{C'}$ is a union of several
pairs of varieties, exchanged by the involution $\iota$.  In general,
some components of $\wt C$ map to themselves). Consider a commutative
diagram:

\xymatrix{
  &&&&X' \ar[r]\ar[d] & X\ar[d] \\
  &&&& \wY \ar[r] &Y }

where $X$ and $X'$ are g.d.c. and $S_2$ varieties, the vertical arrows
are $G$-covers, $X'\to \wY$ is a cover as in the previous section, and
$X'\to X$ is a birational morphism.

 We denote by $B,B'$ the preimages of
$C,C'$ in $X,X'$, and by $\wt{B'}$ the normalization of $B'$.

\begin{equation}\label{big-diagram} 
\xymatrix{
  && \wt{B'} \ar[rr] \ar[dd] \ar@(ul,dl)[]_{j} 
  &&B' \ar[rr] \ar[dd] \ar@{^(->}[dr] 
  &&B \ar@{^(->}[dr] \ar[dd] \\
  && &&&X' \ar[rr] \ar[dd] 
  &&X \ar[dd] \\
  && \wt{C'} \ar[rr] \ar@(ul,dl)[]_{\iota} 
  &&C' \ar[rr] \ar@{^(->}[dr] 
  &&C \ar@{^(->}[dr] \\
  && &&&\wY \ar^{\nu}[rr] 
  &&Y
}
\end{equation}

We first give two constructions for the cover $X\to Y$ starting with
$X'\to\wY$ and the appropriate data for the double locus. One
construction proceeds by $S_2$-fication of the ``nice'' part. The
second one is by a gluing procedure, and the result is very convenient
for computing the invariants of $X$. Finally, we show that indeed
every $X\to Y$ comes from these constructions.

\begin{theorem}\label{thm:glue} Suppose we are given
  \begin{enumerate}
  \item $Y$, $\wY$, $C'$, $(\wt{C'},\iota)$,
  \item a $G$-cover $X'\to \wY$, with $X'$ an $S_2$ and
    g.d.c. variety,
  \end{enumerate}
  Let $B'\to C'$ be the induced cover and let $\wt{B'}\to B'$ be its
  normalization. 
  
  Then $X'$ can be glued to a cover $X\to Y$ with $X$ g.d.c. and $S_2$
  if and only if  it is generically smooth along $B'$,
  and there exists an involution $j\colon \wt{B'}\to \wt{B'}$
  that covers the involution $\iota\colon \wt{C'}\to \wt{C'}$ and commutes
  with the action of $G$ on $\wt{B'}$.  
\end{theorem}

\begin{proof}[Proof by $S_2$-fication]
  Assume that $X$ exists.  
  Then  the map $\wt{B'}\to X$ induces an involution $j$ as required. In addition, if  $X'$ were not generically smooth  along a component  $F$ of $B'$, then $X$ would have generically  at least three branches along the image of $F$. Thus  these two  conditions on $X'$ are  necessary for the existence of~$X$.  

Next we  show that they  are also sufficient.    We start by identifying  the ``bad locus''. It includes the  singular locus of
  $\wY$, the intersection of branch divisors between themselves and
  with $C'$.
   The image of this bad locus in $Y$ has codimension
  $\ge2$. Let $Y_0$ be its complement, and restrict all varieties and
  covers to $Y_0$. 

  The condition that the involution $j$ commutes with the  $G$-action
  implies that for any irreducible component $F$ of $B'$ the subgroup
  $H$ of elements of $G$  that fix $F$  pointwise is the same as the
  supgroup of elements that fix $jF$ pointwise.  
Since $X'$ is generically smooth  along $B'$, one has (cf.   \cite[\S
1]{Pardini_AbelianCovers})    $H=\bZ_n$ for some $n$ and, working \'etale-locally, $H$ acts locally by  $(x, x_2,\dots x_n)\mapsto (\xi x, x_2\dots x_n)$ near $F$ and by $(y, y_2\dots y_n)\mapsto (\xi^a y, y_2,\dots y_n)$  near $jF$  for some primitive
  root $\xi^n=1$ and $(a,n)=1$. Here $y_i=j^*x_i$, $i=2,\dots n$.

  We glue $X'_0$ along $B_0:=\wt{B'}_0/j=B'_0/\iota$ to obtain a
  variety $X_0$ with a finite morphism to $Y_0$. The $G$-action
  extends to $X_0$,  because $j$ commutes with the $G$-action, and is  of the type (smooth) $\times$ (compatible action of curves), where ``compatible'' means that, working \'etale-locally, 
   $\Z_n$ acts on  $xy=0$  by  $x\mapsto \xi x$, $y\mapsto \xi^a y$ 

  Over the double locus we have $\bK[x,y]/(xy)$ and the ring of
  $\bZ_n$-invariants is $\bK[u,v]/(uv)$, where $u=x^n$ and
  $v=y^n$. Thus, $X_0$ has only normal crossings and $X_0\to Y_0$ is a
  $G$-cover. 

  Finally, we apply Lemma~\ref{lem:remove-codim2} to obtain an $S_2$
  and g.d.c. cover $X\to Y$ by taking $S_2$-fication.
\end{proof}

\begin{proof}[Proof by explicit gluing]
  We obtain $X$ by gluing $X'$ along the involution $j\colon\wt{B'}
  \to \wt{B'}$, i.e. as the pushout of the following commutative
  diagram:

  \xymatrix{
    && \wt{B'} \ar[r] \ar[d] &\wt{B'}/j 
    &&\cO_{\wt{B'}}  &\cO_{\wt{B'}/j} \ar[l]
    \\
    && X'
    &&&\cO_{X'} \ar[u]
  }

  Since all varieties are affine over $Y$, $\cO_X$ is the fiber product
  of the corresponding diagram of $\cO_Y$-algebras, in which we
  identify sheaves with their pushforwards on~$Y$.
  We can rewrite this fiber product diagram by saying that $\cO_X$ is the
  kernel in the exact sequence
  \begin{displaymath}
    0\to \cO_X \to \cO_{X'}\oplus \cO_{\wt{B'}/j} 
    \xrightarrow{\beta} \cO_{\wt{B'}}.
  \end{displaymath}
  Further, we have 
  \begin{displaymath}
    0\to \cO_{\wt{B'}/j} \to \cO_{\wt{B'}} \to \cA \to 0,
  \end{displaymath}
  where $\cA$ is the alternating part (if $\chr\bK\ne2$ then
  $\cO_{\wt{B'}}=\cO_{\wt{B'}/j} \oplus \cA$), and the image of
  $\beta$ contains $\cO_{\wt{B'}/j}$. Hence, we have induced exact
  sequences
  \begin{equation}\label{eq:gluing}
    0\to \cO_X \to \cO_{X'}
    \xrightarrow{\alpha} \cA,
    \qquad
    0\to \cO_X \to \cO_{X'}
    \xrightarrow{\alpha} \im\alpha\to 0 
  \end{equation}
  The thus defined variety $X$ is $S_2$ by the next
  Lemma~\ref{lem:sings-after-gluing}, since $\im\alpha$ is a subsheaf of
  $\cA$ and so obviously does not have embedded primes. 
  It is g.d.c. again by looking in
  codimension 1 as in the previous proof. The $G$-action on $X'$
  descends to a $G$-action on $X$ since $j$ commutes with the
  $G$-action on $\wt{B'}$ 
  and  
 by construction the
  subalgebra of $G$-invariants is the algebra of $\wt Y$ glued along
  $\wt{C'}/\iota$, i.e. $\cO_Y$. 
\end{proof}

The varieties $X$ obtained in the two proofs coincide, since they both
have finite morphisms to $Y$, are both $S_2$ and they coincide over an
open subset $Y_0\subset Y$ with $\codim (Y\setminus Y_0)\ge 2$. 

\begin{warning}\label{warn:glue}
  It may happen that there is no covering involution of $B'$ 
  but only
  of its normalization $\wt B'$. Then the double locus of $X$ is
  obtained from $\wt B'/j$ by some additional gluing in codimension 1
  (codimension 2 for $X$). As a consequence, branches of $X$ 
  may not be  $S_2$.
But the variety
  $X$ is $S_2$.  \cite[\S 5.4]{AlexeevPardini_CB} contains multiple
  examples of this phenomenon.
  
  On the other hand, the involution $j$ need not be unique. For
  instance, if $g\in G$ has order 2, then $jg$ is another involution
  satisfying the assumptions for gluing.  The next example shows that
  gluing via different involutions can give rise to non isomorphic
  covers.  
\end{warning}
\begin{example}\label{ex:pinch}

  Let $Y=\{u^2-wv^2=0\}\subset \bA_{u,v,w}$. 
The normalization of $Y$
  is the map $\wt Y=\bA^2_{s,t}\to Y$ defined by $u=st, v=t,
  w=s^2$. Here $C=\{u=v=0\}$, $\wt{C'}=C'=\{t=0\}$ and the involution
  $\iota$ of $\wt{C'}$ is given by $s\mapsto -s$.

Let $X'=\{\epsilon^2=1\}\subset \bA^3_{s,t, \epsilon}$  and let $p\colon X'\to \wt Y$ be the trivial $\Z_2$ cover, given by the projection  on the coordinates $s,t$. The $\Z_2$-action is $\epsilon \mapsto -\epsilon$ and 
 $B'=\wt{B'}=\{t=0, \epsilon^2=1\}$. There are two involutions of  $\wt{B'}$ that lift $\iota$, namely $j_1:=(s,\epsilon)\mapsto (-s, \epsilon)$ and  $j_2:=(s,\epsilon)\mapsto (-s, -\epsilon)$. The cover $X_1\to Y$ obtained by gluing via $j_1$ is obviously the trivial $\Z_2$-cover. 
 
 We describe the cover $X_2\to Y$ obtained by gluing via $j_2$
 following the second proof of Theorem \ref{thm:glue}. The map
 $\wt{B'}\to \wt{B'}/j_2$ corresponds to the inclusion
 $\bK[s\epsilon]\to \bK[s,\epsilon]/(\epsilon^2-1]$ and the map
 $\wt{B'}\to X'$ corresponds to the surjection
 $\bK[s,t,\epsilon]/(\epsilon^2-1) \to
 \bK[s,\epsilon]/(\epsilon^2-1]$. The fiber product of these two ring
 maps can be identified with $R:=\bK[s, t, \epsilon
 t]/(\epsilon^2-1)\subset \bK[s,t,\epsilon]/(\epsilon^2-1)$. The map
 $R\to \bK[x,y,z]/(x^2-y^2)$ defined by $s\mapsto z$, $t\mapsto x$,
 $\epsilon t\mapsto y$ is an isomorphism, hence $X_2$ is the union of
 two copies of $\bA^2$ glued along a line. The cover $X_2\to Y$ is
 given by $(x,y,z)\mapsto (x,yz,z^2)$ and the $\Z_2$-action on $X$ is
 given by $(x,y,z)\mapsto (x,-y,-z)$, thus $(0,0,0)\in Y$ is the only
 branch point.  So the ramification locus of a standard $G$-cover has
 always pure codimension 1 but this not true for the $G$-covers
 obtained from a standard cover by gluing and the analogue of Lemma
 \ref{lem:inertia} does not hold.
\end{example}
  
\begin{lemma}\label{lem:sings-after-gluing}
  With the notations as in the 2nd proof by gluing, assume that $X'$ is
  $S_n$ for some $n\ge2$. Then $X$ is $S_n$ iff $\im\alpha$ is
  $S_{n-1}$. 
\end{lemma}
\begin{proof}
  We use the cohomological interpretation of depth using local
  cohomology \cite[3.8]{Hartshorne_LocalCohomology} (alternatively and
  equivalently one can use $\Ext^i(\cO_{X,Z}/m_{X,Z},\bullet)$). A
  sheaf $\cE$ satisfies $S_n$ iff for every irreducible subvariety
  $Z\subset Y$ one has $H^i_Z(\cE)=0$ for all $i<\min(n,\codim
  Z)$. Looking at the long exact sequence of cohomologies
  corresponding to the short exact sequence~\eqref{eq:gluing}, we get 
  $H^i_Z(\cO_X)=H_Z^{i-1}(\im\alpha)$ for all $i<\min(n,\codim Z)$. The
  statement now follows.
\end{proof}

We spell out Theorem \ref{thm:glue} in a special case, which is of
interest to us because of the applications in \cite{AlexeevPardini_CB}.
\begin{example}\label{ex:2surfaces}
  Take $G=\Z_2^r$. For simplicity of exposition, we assume that
  $Y=Y_1\cup Y_2$ is the g.d.c.  union of two smooth projective
  surfaces that intersect along a smooth rational curve $C$, but all
  our considerations generalize straightforwardly to the case of a
  g.d.c. surface with smooth components whose double locus is a union
  of smooth rational curves.

  We have $\wY=Y_1\sqcup Y_2$, hence an $S_2$ and g.d.c.  $G$-cover
  $X'\to \wY$ is the disjoint union of $S_2$ and g.d.c.  covers
  $\pi_i\colon X'_i\to Y_i$, $i=1,2$. By Corollary
  \ref{cor:structureZ2}, the covers $\pi_i$ are standard. We denote by
  $D^{(i)}_{1},\dots D^{(i)}_{r_i}$, $g^{(i)}_1,\dots g^{(i)}_{r_i}$
  the branch data of $\pi_i$, $i=1,2$. 
  We write $\wt{C'}=C'=C_1'\sqcup C_2'$, $B'=B'_1\sqcup B'_2$ and $\wt{B'}=\wt{B'_1}\sqcup \wt{B'_2}$. We denote by
  $\ga_i$ the generator of 
  subgroup $H_{C'_i}$.   An involution $j$  of $\wt{B'}$ as in
  Theorem \ref{thm:glue} exists if and only if there is an isomorphism
  $\wt{B'_1}\to \wt{B'_2}$ compatible with the $G$-action. This is equivalent to the following
  conditions:
  \begin{enumerate}
  \item $\ga_1=\ga_2=:\ga$,
  \item for $y\in C$, denote by $m^{(1)}_{y,h}$ the intersection
    multiplicity at $y$ of $D^{(1)}_h$ with $C=C_1$, $h=1,\dots r_1$
    and by $m^{(2)}_{y,s}$ the intersection multiplicity at $y$ of
    $D^{(2)}_s$ with $C=C_2$, $s=1,\dots r_2$. Then:

$$\sum_h m^{(1)}_{y, t} g^1_h=\sum_s m^{(2)}_{y, s }g^2_s \quad  \mod \ga, \quad \forall y\in C.$$
\end{enumerate}
Indeed, (1) follows immediately by the fact that $j$ commutes with the
action of $G$. In addition, by the normalization algorithm of \cite[\S
3]{Pardini_AbelianCovers} condition (2) is equivalent to requiring
that the branch data of the normalizations $\wt{B'_1}\to C$ and $\wt{B'_2}\to C$ of the $G/\!\langle\ga\rangle$-coverings of $C=C_1=C_2$ induced by
$\pi_1$ and $\pi_2$ are the same. Since $C$ is smooth rational, the
branch data are enough to determine the building data (cf. Remark
\ref{rem:branchdata}). Since $C$ is projective, the building data
determine the cover up to isomorphism by Proposition \ref{prop:bdata}.

Assume 
that the gluing conditions are satisfied.
Giving an involution of $\wt{B'}$ that commutes with the $G$ action is
the same as giving an isomorphism of $G$-covers $\alpha\colon
\wt{B'_1}\to \wt{B'_2}$.  Then any other such map $\alpha'$ is equal
to $\alpha g$ for some $g\in G$ and the automorphism of $X'=X'_1\sqcup
X'_2$ defined by $x\mapsto x$ if $x\in X'_1$ and $x\mapsto gx$ if
$x\in X'_2$ induces an isomorphism of the cover of $Y$ obtained by
gluing via $\alpha$ with the one obtained by gluing via $\alpha'$. So
in this case all the possible involutions give isomorphic covers.
\end{example}

\begin{theorem}[The reverse]\label{thm:the-reverse}
  Vice versa, every $G$-cover $X\to Y$ with g.d.c. $S_2$ varieties
  $X,Y$ is obtained  via  the gluing construction of Theorem~\ref{thm:glue}.
\end{theorem}

\begin{proof}
  Given $X\to Y$ and the normalization $\wY\to Y$, let $X''$ be the
  fiber product $X''=X\times_Y\wY$. We define $X'$ as
  $X':=S_2(X''\ured)\to X''\ured\to X''$.  Thus, $X'$ is $S_2$ by
  definition, and it maps  to $\wY$. By the universality of
  taking the reduced part and $S_2$-fication, there is an induced $G$-action on $X'$.
  By the universal property of $G$-quotients, we also have a morphism
  $X'/G\to Y$. We claim that it is an isomorphism.

 It is enough to check this in codimension one over
  the double locus. We claim that generically over the double locus of
  $Y$, the cover is (smooth) $\times$ (admissible action of curves),
  where ``admissible'' means that, working \'etale-locally,  $X$ is  given by $xy=0$, and the
  action is $x\mapsto \xi x$, $y\mapsto \xi^a y$ for some primitive
  root $\xi^n=1$ and $(a,n)=1$. Indeed, let $H_F$ be the subgroup of
  elements that restrict to  the identity on  an
  irreducible component $F$ of the double locus of $X$. Then on the
  normalization on both branches we have the same subgroup for the
  preimages $F'$ and $jF'$. Since generically $F',jF'$ are smooth,
  $H_F=\bZ_n$ for some $n\ge1$ (note that one possibly has $n=1$).

  Thus, \'etale locally  the morphism $X\to Y$ can be written as  
  \begin{displaymath}
    \text{(smooth)} \times \bK[u,v]/(uv)\to \bK[x,y]/(xy), 
    \qquad
    u\mapsto x^n,\  v\mapsto y^n, 
  \end{displaymath}
  where  $G$ acts as $x\mapsto \xi x$, $y\mapsto \xi^a y$, $\xi^n=1$,
  $(a,n)=1$. Computing, we get that $X''$ corresponds to  (smooth)$\times \bK[x,y]/(xy,y^n) \oplus \bK[x,y]/(xy,x^n)$, and $X'$ to $\bK[x]\oplus
  \bK[y]$. The quotient $X'/G$ is then $\bK[u]\oplus \bK[v]$, i.e. $\wY$. 

  This proves that $\phi\colon X'/G\to \wY$ is an isomorphism outside a closed
  subset of codimension $\ge2$. Since both are finite over $Y$ and
  $S_2$, $\phi$ is an isomorphism.
\end{proof}

\section{Singularities of covers}\label{sec:geometry}

\subsection{The canonical divisor and slc singularities}
\label{subsec:can-divisor}

Let $Z$ be a variety, let $B_j$, $j=1,\dots n$, be effective Weil
divisors on $X$, possibly reducible, and let $b_j$ be rational numbers
with $0\le b_j\le 1$.
Set $B=\sum_j b_jB_j$.

\begin{definition}\label{defn:lc}
  Assume that $Z$ is a \emph{normal} variety. Then $Z$ has a canonical
  Weil divisor $K_Z$ defined up to linear equivalence. The pair
  $(Z,B)$ is called \emph{log canonical} if
  \begin{enumerate}
  \item $K_Z+B$ is $\mathbb Q$-Cartier, i.e. some positive multiple is
    a Cartier divisor, and
  \item Every prime divisor of $Z$ has multiplicity $\le 1$ in $B$ and
    for every proper birational morphism $h\colon Z'\to Z$ with normal
    $Z'$, in the natural formula
    \begin{displaymath}
      K_{Z'} + h_*\inv B = h^*(K_Z+B) + \sum a_i E_i
    \end{displaymath}
    one has $a_i\ge -1$.
    Here, $E_i$ are the irreducible exceptional divisors of $h$, the
    pullback $h^*$ is defined by extending $\mathbb Q$-linearly the
    pullback on Cartier divisors, $h_*\inv B=\sum b_j h_*\inv B_j$ is
    the strict preimage of $B$.
    The coefficients $a_i$ are called \emph{discrepancies}. For the
    non-exceptional divisors, already appearing on $Z$, one defines
    $a(B_j) = -b_j$.

    If $\chr \bK=0$, then $Z$ has a resolution of singularities
    $h\colon Z'\to Z$ such that $\Supp(h_*\inv B)\cup E_i$ is a normal
    crossing divisor; then it is sufficient to check the condition
    $a_i\ge-1$ for this morphism $h$ only.
  \end{enumerate}
\end{definition}

\begin{definition}
  A pair $(Z,B)$ is called \emph{semi log canonical} if
  \begin{enumerate}
  \item $Z$ satisfies Serre's condition $S_2$,
  \item $Z$ is g.d.c., and no divisor $B_j$ contains any component of
    the double locus of $Z$,
  \item some multiple of the Weil $\mathbb Q$-divisor $K_Z+B$, well
    defined thanks to the previous condition, is Cartier, and
  \item denoting by $\nu\colon \wZ\to Z$ the normalization, the
    pair $(\wZ,\ \text{(double locus)} + \nu_*\inv B )$ is log
    canonical.
  \end{enumerate}
\end{definition}

\begin{lemma}\label{lem:adjunction}
  Let $f\colon X\to Y$ be a finite morphism of degree $d$
  between equidimensional $S_2$ varieties. Assume that either $\chr
  \bK=0$, or $f$ is Galois and $\chr \bK$ does not divide $d$. 

  Let $Y_0$ be an open subset and denote by $f_0\colon X_0\to Y_0$ the
  induced cover. Assume that:
  \begin{itemize}
  \item $\codim(Y\setminus Y_0)\ge2$ and both $X_0$ and $Y_0$ are d.c.,
  \item there exist effective $\bQ$-divisors $B^X$ of $X$ and $B^Y$ of
    $Y$, not containing any component of the double locus, such that
    $(f_0)^*(K_{Y_0}+ B^{Y_0}) = (K_{X_0}+ B^{X_0})$, where $B^{Y_0}$
    is the restriction of $B^Y$ to $Y_0$ and $B^{X_0}$ is the
    restriction of $B^X$ to $X_0$.
  \end{itemize}
  Then:
  \begin{enumerate}
  \item $K_Y+B^Y$ is $\bQ$-Cartier iff so is $K_X+B^X$.
  \item The pair $(Y,B^Y)$ is slc iff so is the pair $(X,B^X)$.
  \end{enumerate}
\end{lemma}
\begin{proof}
  (1) Let $i\colon X_0\to X$ be the inclusion map. If the sheaf
  $L=\cO_Y(N (K_Y+B^Y))$ is invertible then we have a homomorphism
  $$\cO_X(N (K_X+B^X)) = i_*(\cO_{X_0}(N (K_{X_0}+B^{X_0}))) \to f^*L$$
  which is an isomorphism outside of codimension 2. So it must be an
  isomorphism by the $S_2$ condition. Similarly, if the sheaf
  $L'=\cO_X(N (K_X+B^X))$ is invertible then the sheaf $L=\cO_Y(Nd
  (K_Y+B^Y))$ is isomorphic to the norm of $L'$, so is invertible. 

  (2) Assume first that $X$ and $Y$ are normal. In the case this
  statement, due to Shokurov, is very well known. We recall the proof
  because usually it is only stated and proved in characteristic zero.
  Let $h_Y\colon Y'\to Y$ be some partial resolution with normal $Y'$,
  $X'$ be the normalization of $X\times_Y Y'$, and let $h_X\colon
  X'\to X$, $f'\colon X'\to Y'$ be the induced maps.

  Pick an irreducible divisor $E$ on $Y'$, and let $F$ be an
  irreducible divisor on $X'$ over it. By our condition on $\chr \bK$,
  the field extension $\bK(F)/\bK(E)$ is separable, and if $\pi_X$,
  $\pi_Y$ are uniformizing parameters in the DVRs $\cO_{X',F}$ and 
  $\cO_{Y',E}$, then one has $\pi_Y = u \cdot \pi_X^e$ for a unit $u$
  and some integer $e$ dividing $d$ and hence coprime to $\chr \bK$. 

  Then Riemann-Hurwitz formula applies and says that generically along
  $E$ and $F$ one has $(f')^* (K_{Y'}+E) = K_{X'}+F$. Comparing this
  to the identity $(f')^* h_Y^*(K_Y+B^Y) = h_X^* (K_X+B^X)$ and the
  definition of the log discrepancy, one obtains that $1+a_F =
  e(1+a_E)$. Thus, $a_F\ge -1$ $\iff$ $a_E\ge -1$. This proves that
  $(X, B^X)$ is lc iff $(Y,B^Y)$ is lc.

  Now consider the general g.d.c. case. Let $\nu_X\colon \wX \to
  X$ be the normalization. We have
  \begin{displaymath}
    K_{\wX} + B^{\wX} :=
    \nu_X^*(K_X + B^X) = K_{\wX} + \nu_{X*}\inv B^X 
    + {\rm (double\ locus)}, 
  \end{displaymath}
  and similarly for $Y$. Thus, the double loci appear in the divisors
  $B^{\wX}$, $B^{\wY}$ with coefficient 1. By the Riemann-Hurwitz
  formula again, for the normalizations we still have $\tilde f^*
  (K_{\wY} + B^{\wY}) = K_{\wX} + B^{\wX}$. We conclude by applying
  the normal case.
\end{proof}

We now extend Definition \ref{defn:hurwitz} of Hurwitz divisor to the
case of a g.d.c.  base $Y$:

\begin{definition}\label{defn:hurwitz2}
   Let $\pi\colon X\to Y$ be a $G$-cover of $S_2$ and
   g.d.c. varieties. For a prime Weil divisor $F\subset Y$, we define $\rho_F\in \Q$ as follows:
   \begin{itemize} 
   \item if $F$ is contained in the double locus of $Y$, then $\rho_F=0$;
   \item if $F$ is not contained in the double locus of $Y$, but $\pi\inv(F)$ is contained in the double locus of $X$, then $\rho_F=1$,
   \item  if $F$ is not contained in the double locus of $Y$,
     $\pi\inv(F)$ is not 
contained in the double locus of $X$ and $m$ is the ramification order of $\pi$ at $F$, then $\rho_F=\frac{m-1}{m}$.
     \end{itemize}
      We define the Hurwitz divisor $D$ of $\pi$  to be the $\Q$-divisor  $\sum_{F} \rho_FF$.
      
      Notice that if $X\to Y$ is a standard $G$-cover with $X$ g.d.c. this definition coincides with Definition \ref{defn:hurwitz} by Theorem \ref{thm:structure}.
\end{definition}

Note that $D$ does not contain any components of the double locus of
$Y$.

\begin{proposition}\label{prop:interesting-case}
  Let $\pi\colon X\to Y$ be a $G$-cover as in Definition
  \ref{defn:hurwitz2} and let $D$ be the Hurwitz divisor of $\pi$, let
  $X'\to \wY$ be the corresponding $S_2$ and g.d.c. $G$-cover (cf.\S~
  \ref{ssec:Ynonnormal}) 
  Then
  \begin{enumerate}
  \item $K_X$ is $\bQ$-Cartier iff so is $K_Y+D$, and then
    $K_X=\pi^*(K_Y+D)$. 
  \item $X$ is slc iff so is the pair $(Y,D)$.
  \end{enumerate}
\end{proposition}
\begin{proof} 
Recall that $|G|$ and $\chr\bK$ are coprime by assumption.
So  Lemma \ref{lem:adjunction} applies and we may
  assume that $Y$ is d.c..  We need to show that
$K_{X}=\pi^*(K_Y+D).$
This is equivalent to the following equality for the cover $\tilde\pi\colon
\wX \to \wY$, where $\wX$ is the normalization of $X'$ (and of $X$):
\[ K_{\wX}+(\rm{double\ locus})={\tilde\pi}^*(K_{\wY}+(\rm{double\
  locus})+\nu^*D).\] 
In view of Definition \ref{defn:hurwitz2} the formula follows easily by
the usual Hurwitz formula.
\end{proof}
\subsection{Cohen--Macaulay covers}

By Lemma~\ref{lem:flatness-and-CM}, a $G$-cover over a smooth base is
CM. Here, we give a partial generalization of this case to the case of
a non-normal base. We use the notations of Theorem~\ref{thm:glue} and the
exact sequence \eqref{eq:gluing}.

\begin{prop}\label{prop:CM}
  Assume that $X'$ is CM (for example, $\wY$ is smooth). Then $X$ is
  CM iff the sheaf $\im\alpha$ is CM.
\end{prop}
\begin{proof}
  Immediate by Lemma~\ref{lem:sings-after-gluing}.
\end{proof}
Using Proposition \ref{prop:CM} it is not hard to give examples of
abelian covers $X\to Y$ such that $Y$ is CM and g.d.c., and $X$ is
g.d.c. and $S_2$ but not CM:
\begin{example}\label{ex:nonCM}
We take $G=\Z_2$ and assume $\chr \bK\ne 2$;  for  any prime $p$ one can construct similar examples  with $G=\Z_p$ and $\chr \bK\ne p$. 

Let $Y=Y_1\cup Y_2$ be the union of 2 copies of $\pp^3$ glued transversally  along a plane $C$. Let $L_1$ and $L_2$ be distinct lines on $C$  and for $i=1,2$ let $D_i\subset Y_i$ be a  quadric that restricts to $2L_i$ on $C$.  For a generic choice, $D_i$ is a quadric cone with vertex $y_i\in L_i$  and the points $y_1$, $y_2$ and $y_3:=L_1\cap L_2$ are distinct. Let $X'_i\to Y_i$ be the double cover of $Y_i$ branched on $D_i$ and let $X'=X'_1\sqcup X'_2$. Then $X'$ is Gorenstein,   has an ordinary double point over $y_1$ and $y_2$  and no other singularity. 
Write $C'=C'_1\sqcup C'_2$ and $B'=B'_1\sqcup B'_2$; then  $B'_i$ is the union of two copies of $C'_i$ glued transversally  along $L_i$ and $\wt{B'}\to C'$ is the trivial $\Z_2$-cover. Hence  there exists an involution $j$ of $\wt{B'}$ that commutes with the $\Z_2$-action, and by Theorem \ref{thm:glue} $X'$ can be glued to an $S_2$ and g.d.c. cover $X\to Y$.
The  d.c. locus of  $X$  is the complement of the preimage  of $L_1\cup L_2$.

In the exact sequence  \eqref{eq:gluing} each term splits under the $G$-action and the maps are compatible with the splitting,  so we get two exact sequences, one for each character of $G$.  
Since $\cA=\OO_C\oplus \OO_C$ and $\Z_2$ acts on $\cA$ by switching the two summands,  the sequence for the nontrivial character is: 
$$0\to \cF_{-}\to  \OO_{Y_1}(-1)\oplus\OO_{Y_2}(-1) \xrightarrow{\alpha^-} \OO_{C},$$
where $\cF_{-}$ (resp. $\cA^-$) is the antiinvariant summand of $\OO_{X}$ (resp. of $\cA$).
By definition, the map $\OO_{Y_i}(-1)\to \OO_{C}$  factorizes  as $\OO_{Y_i}(-1)\to \OO_{C}(-L_i)\to \OO_{C}$.
Hence, $\im \alpha^{-}$ coincides with $\cI_{y_3}\OO_C$, the maximal ideal of $y_3$ in $C$, and therefore it is not $S_2$. It follows by Proposition \ref{prop:CM} that $X$  is not CM  over $y_3$.

Let $\bar y\in L_1$ be a point distinct from $y_3$; in a neighbourhood
of $\bar y$ we have $(D_1+D_2)\cap Y_2=L_1$, thus $D_1+D_2$ is not
$\bQ$-Cartier. Since $Y$ is Gorenstein, it follows that $2K_Y+D_1+D_2$
is not $\bQ$-Cartier either, hence $K_X$ is not $\bQ$-Cartier by
Proposition \ref{prop:interesting-case}. 
\end{example}
\subsection{Cartier index of $K_X$}\label{sec:gorenstein}

All the statements in this section are \'etale local, 
so we often pass to a smaller neighbourhood of a
point without explicit mention of the fact.

For convenience, we write ``$K_X$'' to denote the divisorial
sheaf $\omega_X$ (recall that $X$ is Gorenstein in codimension 1 and
$S_2$). We also use the additive notation $D_1+D_2$ for the sheaf
$\left( \cO_X(D_1) \otimes \cO_X(D_1) \right)^{**}$.

\subsubsection{Standard covers with $Y$ normal}
 \label{subsec:gor-standard}

We consider the following situation:
\begin{itemize}
\item $Y$ is a  normal variety  and $C$ is a reduced effective  divisor on $Y$ such that $K_Y+C$ is Cartier;
\item $\pi\colon X\to Y$ is a standard g.d.c.  $G$-cover (so $X$ is automatically $S_2$ by Lemma \ref{lem:flatness-and-CM}). We  assume that $X$ is generically smooth over  $C$ and we denote by $B$ the preimage of $C$ in $X$. So $B$ is also a reduced effective divisor. 
\end{itemize}

Let $D$ be the Hurwitz divisor of $\pi$; then we have:
$$K_X+B=\pi^*(K_Y+D+C).$$
Thus,  if $d$ is  the exponent
of $G$,  then the divisor $d(K_Y+D+C)$ is Cartier 
(recall that the divisors $D_i$ are Cartier by the definition
of a standard cover in Section~\ref{ssec:standard}) and thus 
$d(K_X+B)$ is also Cartier.  

Fix a point $y\in Y$;  the purpose of this section is to compute the
Cartier index of $K_X+B$ at a point $x\in X$ such that $\pi(x)=y$. 
 Here we are interested mainly in the case $B=0$, but the 
 case of a pair
is needed in the next section to treat the case  $Y$  non-normal.

In order to state our result we need some notation.  We label the
branch data $D_i, (H_i,\psi_i)$, $i=1,\dots k$, in such a way that
$D_i\subseteq C$ iff $i\le p$.
Since the question is local on $Y$ we may assume that $y\in D_i$ for
every $i$. Consider the map $\ol G:=\oplus H_i\to G$. By Lemma
\ref{lem:inertia} the image of this map is the inertia subgroup $H_y$;
we denote by $N$ the kernel. We let $\ol {\chi}\in {\ol G}^*$ be the
character $\psi_{p+1}\cdots \psi_k$.

\noindent\underline{Reminder:} Since the group $G$ is finite abelian,  the map $G^*\to H_y^*$ is surjective. So the character $\ol{\chi}$ is the pullback of a character of $H_y$ iff it is the pullback of a character of $G$.

\begin{proposition}\label{prop:index-normal}
Notation and assumptions as above. \\
The Cartier index of $K_X+B$ at $x$ is equal to the order of $N/(N\cap
\ker\ol{\chi})$.

In particular, $K_X+B$ is Cartier iff $\ol{\chi}$ is the pull back of a character $\chi\in G^*$.
\end{proposition}
\begin{proof}
Since the question is  local, we  may assume that the line bundles $L_{\chi}$,
$\OO_Y(D_i)$ and $\OO_Y(K_Y+C)$ 
 are trivial. The map $X\to X/H_y$ is \'etale, hence up to replacing $Y$ by $X/H_y$ we may assume that $H_y=G$, or, equivalently, that $\pi\inv(y)=\{x\}$.  We denote by $u_1, \dots u_k$
local equations of $D_1,\dots D_k$ near $y$. By Remark
\ref{rem:local}, up to passing to an \'etale cover of $Y$   we may assume
that $X$ is given by:
\begin{equation}
  z_{\chi}z_{\chi'}=\Pi_1^k u_i^{\epsi^i_{\chi,\chi'}} z_{\chi\chi'}, \quad \chi,\chi'\in G^*\setminus\{1\},
\end{equation}
The equations:
\begin{equation}\label{eq:X0}
  z_1^{m_1}=u_1,\quad \dots \quad z_k^{m_k}=u_k
\end{equation}
define inside $Y\times \bK^k$   a  $\ol G$-cover $\ol X\to Y$ ($\ol G$
acts on $z_i$ via the character $\psi_i$), the {\em maximal totally ramified cover} of $Y$ with branch data $D_i,(H_i,\psi_i)$ (here we regard $H_i$ as a subgroup of $\ol G$). Since $Y$ is g.d.c. by assumption and $X\to Y$ and $\ol X\to Y$ have the same  Hurwitz divisor, $\ol X$ is also g.d.c. by Theorem \ref{thm:structure}.

For every $\chi\in G^*$,
write $\chi =\psi_1^{a_{\chi}^1}\cdots \psi_k^{a_{\chi}^k}$, with
$0\le a^{\chi}_i<m_i$ for $i=1,\dots k$; then setting
$z_{\chi}=z_1^{a_{\chi}^1}\cdots z_k^{a_{\chi}^k}$ defines a map
$p\colon \ol X\to X$ which is the quotient map for the action of the kernel 
$N$ of $\ol G\to G$. The map $p$ is unramified in codimension 1 and $p\inv(x)$
consists of just one point $\ol x$.

Denote by $\ol B$ the preimage of $C$ (and of $B$) in $\ol X$; observe that  $K_Y+D+C$ pulls back to $K_X+B$ on $X$ and to  $K_{\ol X}+{\ol B}$ on $\ol X$.
If $\tau$ is a generator of $\OO_Y(K_Y+C)$ then $\OO_{\ol X}(K_{\ol X}+\ol B)$ is generated by the residue $\sigma$ on $\ol X$ of the rational differential form:
$$\frac{(z_1^{m_1-1}\cdots z_p^{m_p-1})dz_1\wedge\dots \wedge dz_k\wedge \tau}{(z_1^{m_1}-u_1)\cdots (z_k^{m_k}-u_k)}.$$
Thus $\OO_{\ol X}(K_{\ol X}+\ol B)$ is invertible and $G$ acts on the local generator $\si$ via the character $\ol\chi$.
Set $Z:=\ol X/(N\cap  \ker \ol \chi)$. The map $\ol X\to Z$ is unramified in codimension 1 and   $\si$ descends on $Z$  to a generator of $\OO_Z(K_Z+B_Z)$, where $B_Z$ is the image of $\ol B$. The map $Z\to X$ is  a cyclic cover with Galois group $N/(N\cap \ker \ol \chi)$ with the following properties:
\begin{itemize}
\item it is unramified in codimension 1 and  the preimage of $x$ consists only of one point,
\item  the pull back of $\OO_X(K_X+B)$ is a line bundle on which the Galois group acts via a primitive character. 
\end{itemize}
It follows that $Z\to X$ is a canonical cover and that  the Cartier index of $K_X+B$ at $x$ is equal to $[N:N\cap\ker \ol\chi]$.
\end{proof}

\begin{corollary}\label{cor:gor-smooth} Let $\pi\colon X \to Y$ be
  a standard abelian with $X$  and $Y$ g.d.c.  and $Y$ 
  Gorenstein, let $y\in Y$ and let $x\in X$ be a point such that $\pi(x)=y$.  Then $X$ is Gorenstein at $x$ iff the character $\ol \chi$
  descends to a character $\chi$ of $H_y$.
\end{corollary}
\begin{proof} $X$ is Cohen-Macaulay by Lemma \ref{lem:flatness-and-CM} and $K_X$ is Cartier by 
  Proposition \ref{prop:index-normal}.
   \end{proof}
\begin{remark}
  Corollary \ref{cor:gor-smooth} is proven in \cite{Iacono} under the
  assumption that $X$ is  normal and $Y$ is smooth.
\end{remark}

\subsubsection{The case $Y$ non-normal}\label{ssec:gor_nonnormal}
 
 Here we consider the problem of determining the Cartier index of $K_X$ at a point $x\in X$  of a $G$-cover $X\to Y$ with $Y$ non-normal of Cartier index 1. The situation is much more complicated than in the case $Y$ normal  and we are able to give only a partial answer, that is however sufficient for the applications in  \cite{AlexeevPardini_CB}. The main difficulty is that one does not know how to write down an analogue of the maximal totally ramified cover used in the proof of Proposition \ref{prop:index-normal}.

 \smallskip 
 
 We consider the following setup:
 \begin{itemize}
 \item $Y=Y_1\cup \dots \cup Y_t$, where  $Y_i$ is irreducible for $i=1,\dots t$,  is a g.d.c. and $S_2$ variety; $\wY= \wt{Y}_1\sqcup\dots \sqcup \wt{Y_t} \to Y$ is the normalization,
 \item $\pi\colon X\to Y$ is  an $S_2$ and g.d.c. $G$-cover obtained by
gluing a cover $X'=X'_1\sqcup\dots \sqcup X'_t\to \wY$ such that $X'_i\to \wt Y_i$ is standard for every $i$,  
\item $y\in Y$  and  $x\in X$ are points such that $\pi(x)=y$; we assume that $y$ lies on every component of the branch locus of $\pi$. \end{itemize}

 We denote by $D_i, (H_i,\psi_i)$, $i=1,\dots k$ the branch data of the standard cover $X'\to\wt Y$ and we assume that $D_i$ is contained in the preimage $C'$ of the double locus of $Y$ if and only if $i\le p$. Consider the map $\ol G:=\oplus H_i\to G$.    
  As in the case $Y$ normal, we denote by $\ol \chi\in {\ol G}^*$ the character $\psi_{p+1}\cdots \psi_k$. 
  Then:
  \begin{proposition} \label{prop:index-non-normal}
  In the above setup, if $K_X$ is Cartier, then:
  \begin{enumerate}
  \item $K_Y+D$ is $\Q$-Cartier,
  \item  $\ol \chi$ is the pullback of a character $\chi\in G^*$.
  \end{enumerate}
  \end{proposition}
  \begin{proof} (1) Follows immediately by Proposition \ref{prop:interesting-case}.
  
  (2) For every $i=1,\dots t$ denote by $C'_i\subset \wt{Y_i}$
  (resp. $B'_i\subset X'_i$) the preimage of the double locus of $Y$ in $\wt{Y_i}$ (resp. in $X'_i$).
  Let $\chi\in G^*$ be the character via which $G$ acts on $\OO_X(K_X)\otimes\bK(x)$   at $x$. Let  $x'_i\in X'_i$ be  a point that maps to $x$ and let $y_i$ be the image of $x'_i$ in $\wt{Y_i}$. Since  $K_X$  pulls back to $K_{X'_i}+B'_i$ on $X'_i$, the inertia subgroup $H_{y_i}$ acts on  $\OO_{X'_i}(K_{X'_i}+B'_i)\otimes\bK(x'_i)$  via the restriction of $\chi$. 
  Set $\ol G_{y_i}:=\oplus_{\{j|y_i\in  D_j\}}H_j$ and let $\ol \chi_{y_i}$ be  the restriction of $\ol \chi$ to $\ol G_{y_i}$; the map $\ol G_{y_i}\to H_{y_i}$ is a surjection by Lemma \ref{lem:inertia}.  By the proof of   Proposition \ref{prop:index-normal}  $\chi$ pulls back on  $\ol G_{y_i}$ to $\ol \chi_{y_i}$. Since $\ol G=\sum_{\{y'\in \wt Y| y'\mapsto y\}}  \ol {G}_{y'}$, it follows that $\chi$ pulls back to $\ol \chi$ on $\ol G$.
  \end{proof}

 We now
prove a partial converse of Proposition \ref{prop:index-non-normal}.
 Assume that for every component $\wt{Y_i}$ 
of $\wY$ the map $\wY \to
Y$ induces a homeomorphism $\wt{Y_i}\to Y_i$ onto its image
(this is always true up to an \'etale cover). 
Then we
associate to $(Y,y)$ an incidence graph $\Ga_{Y,y}$ as follows:
 
-- the vertices of $\Ga_{Y,y}$ are indexed by the branches of $(Y,y)$,

-- the edges are indexed by the components of the double locus $C$ of $Y$,

-- the edge corresponding to a component $F$ of $C$
connects the vertices corresponding to the two branches of $Y$ through
$F$.

\begin{prop}\label{prop:index-graph}
  In the above setup, assume that:
  \begin{enumerate}
  \item the graph $\Ga_{Y,y}$ is a tree,
  \item  $K_Y$ is Cartier and there exists $m$ such that $m(K_Y+D)$ is Cartier and $(m,\chr
    \bK)=1$,
  \item $\ol \chi$ is the pullback of a character $\chi\in G^*$.
  \end{enumerate}
  Then $K_X$ is Cartier.
\end{prop}
\begin{proof}
   Let $C'_i\subset \wt Y_i$ the restriction of the double locus $C'$ of $\wt Y$ and let $B'_i\subset X'_i$ be the preimage of $C'_i$. Let $y_i\in \wt Y_i$ be the only point that maps to $y\in Y$; let $\ol{G}_{y_i}$ and $\ol{\chi}_i$ be defined as in the proof of Proposition \ref{prop:index-non-normal}.
 
   By assumption (3),  the divisor
  $K_{X'_i}+B'_i$ is Cartier by  Proposition \ref{prop:index-normal}.
  By the following Lemma \ref{lem:character},
 up to replacing $(Y,y)$ by an \'etale
  neighbourhood we may assume that for $i=1,\dots t$ the sheaf
  $\OO_{X'_i}(K_{X'_i}+B'_i)$ is trivial and has a generator
  $\si_i$ on which $G$ acts via $\chi$.  By Proposition
  \ref{prop:interesting-case}, there exists a local generator $\tau$
  of $\OO_X(mK_X)$ near $x$. For every $i$, by Lemma
  \ref{lem:character},  $\tau$ pulls back on $X'_i$ to $h_i\si_i^m$  where
  $h_i$ is a nowhere vanishing regular function on $\wt{Y_i}$. Up to  passing to an
  \'etale cover of $Y$  we may assume that $h_i$ has an $m$-th root
  $f_i$ for every $i$. So we may replace $\si_i$ by $f_i\si_i$ and
  assume that $\tau$ pulls back  to $\si_i^m$ for every $i$.

  Now let $U\subset X$ be an open set such that $U$ is d.c. and the complement of $U$ has codimension $>1$.  
  Let $F$ be an irreducible
  component of the double locus $C$  of $Y$ and let $Y_a$, $Y_b$ be the components of $Y$  that contain  $F$. 
    Choose an irreducible component $E$ of
  the inverse image of $F$ in $U$. It makes sense to compare $\si_a$
  and $\si_b$ along $E$, since they both restrict to local generators
  of $\OO_E(K_E)$. Since $\si_a^m=\si_b^m$, there exists $\zeta\in
  \mu_m$ such that $\si_a=\zeta \si_b$ along $E$. Since $G$ acts on
  $\si_a$ and $\si_b$ via the same character $\chi$ and $G$ acts
  transitively on the components of the preimage of $F$,
  $\zeta_F:=\zeta$ depends only on $F$. So $\{\zeta_F\}$ represents
  a class in $H^1(\Ga_{Y,y},\mu_m)$. Since $\Ga_{Y,y}$ is a tree, we
  can find $\lambda_i\in \mu_m$ such that the local generators
  $\lambda_i\si_i$ glue to give a local generator $\si$ of
  $\OO_X(K_X)$ on which $G$ acts via $\chi$.
\end{proof}

We complete the proof of Proposition \ref{prop:index-graph} by proving the following:
\begin{lemma}\label{lem:character}
 Let $Z\to W$ be a standard $G$-cover with building data $L_{\chi}, D_i, (H_i,\psi_i)$.  
 
 Let $w \in W$ be a point and let $H$ be the inertia subgroup of $w$.   Let $L$ be a $G$-linearized line bundle of $Z$, let $z\in Z$ be a point that maps to $w$ and let $\phi\in H^*$ be the character via which $H$ acts  on $L\otimes \bK(z)$.  Then:
\begin{enumerate}
\item  let $\chi\in G^*$ be such that $\chi|_H=\phi$; then, up to replacing  $W$ by an \'etale neighbourhood of $w$,  there exists a
generator $\sigma$ of $L$ such that $G$ acts on $\sigma$ via the character $\chi$;
\item $\sigma$ is uniquely  determined by $\chi$ up to multiplication by a nowhere vanishing regular function of $W$.
\end{enumerate}
\end{lemma}
\begin{proof}
(2) Assume that $\si,\si'$ are generators of $L$ on which $G$ acts via the character $\chi$.  Then $f:=\si/\si'$ is a  regular $H$-invariant function on $Z$,   so it is a function on $W$.
\smallskip

(1) We break the proof into three steps.  

\noindent\underline{Step 1:} {\em the case $H=G$.}\newline
  Let $s$ be a generator of $L$ near $z$. The group $H$ acts on the
  vector space $V$ of local sections of $L$ spanned by the elements
  $h_*s$, $h\in H$. $V$ is finite-dimensional, and decomposes under
  the $G$-action as a direct
  sum of eigenspaces.  Since $s(z)\ne 0$ and $s\in V$, there exists an
  eigenvector $\si\in V$ such that $\si(z)\ne 0$. Since $G$  acts on $L\otimes \bK(z)$ via $\chi$, $\si$ belongs to the eigenspace corresponding to $\chi$. 
 \smallskip

 \noindent\underline{Step 2:} {\em the case in which  $G=H\oplus N$ for some $N$.}\newline
   Consider the factorization $Z\to Z':=Z/N\to W$. The map $Z'\to W$ is an $H$-cover such that the preimage of $w$ consists of one point $z'\in Z'$. 
The subgroup $N$ acts freely on $Z$, hence
   $L$ descends to an $H$-linearized
  line bundle $L'$ on $Z'$.  Then by Step 1 there
  exists a local generator $\si'$ of $L'$ near $z'$ such that
  $H$ acts on $\si'$ via $\phi$. Pulling back to $Z$ we get a
  generator $\tau$ of $L$ on which $H$ acts via $\phi$ and $N$ acts
  trivially.
 
  Denote by $\phi'$ the restriction of $\chi$ to $N$, so that
  $\chi=(\phi, \phi')$.  Consider the factorization $Z\to Z'':=Z/H\to
  W$.  The map  $Z''\to
  W$ is a \'etale $N$-cover. So there exists a nowhere vanishing
  function $f$ on $Z''$ such that $N$ acts on $f$ via the character
  $\phi$. Thus $G$ acts on $\si:=f\tau$ via the character $\chi$.
\smallskip

 \noindent\underline{Step 3:} {\em the general case.}\newline   Choose a
  finite abelian group $N$ with a surjective map $G_0:=H\oplus N\to
  G$ that extends the inclusion $H\to G$ and let $T$ be the kernel
  of $G_0\to G$. By Proposition
  \ref{prop:bdata}, up to replacing $W$ by an \'etale
  neighbourhood of $w$, we may also assume (cf. \eqref{eq:local}) that $Z\to
  W$ is given  inside $W\times \bK^k$ by the equations:
  \begin{equation}\label{eq:localc2}
   y_{\chi}y_{\chi'}=\Pi_1^k u_i^{\epsi^i_{\chi,\chi'}} y_{\chi\chi'}, \quad \chi,\chi'\in G^*\setminus\{1\},
  \end{equation}
  where $u_i$ is a local equation for $D_i$, $i=1,\dots k$.  The
  branch data for $Z$ can be interpreted in an obvious way as branch
  data for a $G_0$-cover.  Letting $Z_0\to W$ be the $G_0$-cover given
  by the equations analogous to \eqref{eq:localc2}, we have $Z= Z_0/T$
  by construction.  Let $L_0$ be the pull back of $L$ to $Z_0$. $L_0$
  has a natural $G_0$-linearization and $H$ is a direct summand of
  $G_0$, hence by Step 2 there exists a generator
  $\si_0$ of $L_0$ on which $G_0$ acts via the character $\chi_0$ of
  $G_0$ induced by $\chi$. Since $T$ acts freely on $Z_0$ and
  $T\subset \ker \chi_0$ by construction, $\si_0$ descends to a
  generator $\si$ of $L$ on  $Z$ on which $G$ acts via $\chi$.
\end{proof}

\section{Slc $\Z_2^r$-covers of surfaces}
\label{sec:Z2-covers}

\subsection{Setup}\label{ssec:setup-Z2}
In this section we make a detailed study of $\Z_2^r$-covers of
surfaces. We use freely the notation introduced in
\S\ref{ssec:Ynormal}. In particular, we refer the reader to the
commutative diagram~(\ref{big-diagram}) and Theorem~\ref{thm:glue}.

 The situation that we consider is the following:
\begin{itemize}
\item $Y$ is a g.d.c.
  surface with smooth irreducible components $Y_1,\dots Y_t$. The
  irreducible components $F_1,\dots F_s$ of the double curve $C$ of
  $Y$ are smooth, $Y$ is d.c. at the smooth points of $C$ and it is
  analytically isomorphic to the cone over a cycle of rational curves
  at the singular points of $C$. In particular, $Y$ is Gorenstein.
 \item $G=\Z_2^r$ and $\pi\colon X\to Y$ is a $G$-cover with $X$ g.d.c. and $S_2$,  obtained as in
  Theorem \ref{thm:glue}
  by gluing a cover $X'\to \wY=Y_1\sqcup \dots \sqcup Y_t$ such that for every $i=1, \dots t$ the restricted  cover $\pi_i\colon X'_i\to Y_i$ is standard  with  building data  $L_{i,\chi}$, $D_{i,j_i}$.

\item The $D_{i,{j_i}}$ and the components of the
   double curve $C'$ are ``lines'' of $Y$, namely they are
  smooth and meet pairwise transversally.
\item The intersection points of the support of the Hurwitz divisor $D$ of $\pi$ with the double
  curve $C$ of $Y$ are smooth points of $C$.
\item  $K_Y+D$ (or, equivalently, $D$, since $Y$ is Gorenstein) is $2$-Cartier and the  pair $(Y,D)$
  is slc, so that by Proposition \ref{prop:interesting-case} $X$ is
  slc and $K_X$ is $2$-Cartier. Recall that, since we assume that the
  components of $\nu^*D$ and of $C'$ are lines, the pair $(Y,D)$ is slc iff on $\wY$
  the divisor $\nu^*D+C$ has components of
  multiplicity $\le 1$ and has multiplicity $\le 2$ at every point.

\item For every $y\in Y$ that is singular for $C$, label the
  components $Y_1,\dots Y_q$ of $Y$ containing $y$ in such a way that
  for every $i=1,\dots q$ the surfaces $Y_i$ and $Y_{i+1}$ meet along
  an irreducible curve $F_i$ containing $y$ (the indices are taken
  modulo $q$) and let $g_i\in G$ be the generator of the inertia
  subgroup of $F_i$. By Theorem \ref{thm:glue}, for every $i$ we have
  $g_{i-1}=g_{i+1} \mod g_i$.  We assume that the natural map $\langle
  g_i\rangle\oplus \langle g_{i+1}\rangle\lra H_y$ is an isomorphism
  for every $i=1,\dots q$.

  These conditions imply that the fibre  of $X\to Y$ over $y$
  consists of $2^r/|H_y|$ points. At each of
  these points $X$ is analytically isomorphic to the cone over a cycle of $q$
  smooth rational curves.
\end{itemize}
 All the above assumptions are satisfied in the cases considered in \cite{AlexeevPardini_CB}.

\subsection{Numerical invariants}

Here we assume that the surface $Y$ is projective.

By Proposition \ref{prop:interesting-case}, $K^2_X$ can be computed as
follows:
\begin{equation}\label{eq:canonical}
  K^2_X=2^r(K_{\wY}+\nu^*D+({\rm double\ locus}))^2=\sum_i 2^r(K_{Y_i}+ D|_{Y_i}+({\rm double\ locus})|_{Y_i})^2.
\end{equation}

 To compute the cohomology of $\OO_X$,  we are going to  write down explicitly in the above  situation the sequences \eqref{eq:gluing} in the second proof of Theorem \ref{thm:glue} (as usual we push forward  to $Y$ all the sheaves). Since all  the maps   are $G$-equivariant, the sequences    \eqref{eq:gluing} split as  sums  of  exact sequences:
\begin{equation}\label{eq:glue-surf}
    0\to \cF_{\chi} \to \oplus_{i=1}^t L_{i,\chi}\inv 
    \xrightarrow{\alpha} {\cA}_{\chi},
    \qquad
    0\to \cF_{\chi} \to \oplus_{i=1}^t L_{i,\chi}\inv 
    \xrightarrow{\alpha} (\im\alpha)_{\chi}\to 0, 
  \end{equation}
where $\chi$ varies  in $G^*$ and $G$ acts in $\cF_{\chi}$, $\cA_{\chi}$ and $(\im\alpha)_{\chi}$ via $\chi$. 

To describe the sheaves $\cA_{\chi}$ and $(\im\alpha)_{\chi}$, we need to introduce  some more notation.
  Given a
component $F_l$ of $C$ we denote by $g_l\in G$ the generator of the
inertia subgroup of $F_l$ and by $Y_{a_l}$ and $Y_{b_l}$ the two
components of $Y$ that contain $F_l$. We denote by $E_l$ (resp. $E_{l, {a_l}}$, $E_{l,{b_l}}$)  the preimages of $F_l$ in $X$ (resp. $X'_{a_l}, X'_{b_l}$)   and by $\wt{E_l}$ the common normalization of  $E_l$, $E_{l, {a_l}}$, $E_{l,{b_l}}$ (cf. Example \ref{ex:2surfaces}). In the following  commutative diagram:

\xymatrix{
  &&&& & \wt{E_l} \ar[dl] \ar[dr] \ar[d]
  & \\
  &&&& E_{l, {a_l}} \ar[dr]\ar[r]  & E_l \ar[d] & E_{l, {b_l}} \ar[dl]\ar[l]\\
  &&&& 
  & F_l
  &
} 
\noindent the  maps to $F_l$ are 
  $G/\!\langle g_l\rangle$-covers and the remaining maps are finite and birational. 
  The cover $E_{l, {a_l}}\to F_l$   is  standard and its  building data can be recovered  from  those  of $X'_{a_l} \to Y_{a_l}$ as follows:
  \begin{itemize}
  \item we identify $(G/\!\langle g_l\rangle)^*$ with $\langle g_l\rangle^{\perp}\subseteq G^*$ and for every $\chi\in \langle g_l\rangle^{\perp}$ we restrict $L_{\chi}^{a_l}$ to $F_l$,
  \item  for every   $D^{a_l}_j$ with $g_j\ne g_l$, we label  each point of $D^{a_l}_j|_{F_l}$ with the image of $g_j$ in $G/\!\langle g_l\rangle$.
  \end{itemize}
  The same can be done of course for $E_{l, {b_l}}\to F_l$.  Let $y\in
  F_l$ be a point such that $\nu^*D$ has multiplicity 1 at the points
  of $\wt Y$ that map to $y$ (since we assume that $2D$ is Cartier,
  the multiplicity of $\nu^*D$ is the same at all points lying over
  $y$).  Recall that by assumption $Y$ is d.c. at $y$; denote by
  $\al_{y, 1}$ $\al_{y,2}$ the elements of $G$ associated to the two
  branch lines of $X'_{a_l}\to Y_{a_l}$ containing $y$ and by $\be_{y,
    1}$, $\be_{y,2}$ the elements of $G$ associated to the two branch
  lines of $X'_{b_l}\to Y_{b_l}$ containing $y$.  We have $\al_{y,
    1}+\al_{y,2}=\be_{y, 1}+\be_{y,2}$ modulo $g_l$ (cf. Example
  \ref{ex:2surfaces}). Then $E_{l, a_l}$ is singular over $y$ iff
  $\al_{y, 1}$ and $\al_{y,2}$ are both different from $g_l$, namely
  iff there exists a character $\chi$ with $\chi(g_l)=1$ and
  $\chi(\al_{1,y})=\chi(\al_{2,y})=-1$.  For each $\chi\in G^*$ and
  $l$ such that $\chi(g_l)=1$ we denote by $A_{l, {\chi}}$ the set of
  points $y\in F_l$ such that
  $\chi(\al_{1,y})=\chi(\al_{2,y})=-1$, and  we take $A_{l, {\chi}}$ to be the empty set if $\chi(g_l)\ne 1$.
  We define in a
  similar way $B_{l,{\chi}}$ by considering the cover $E_{l,{b_l}}\to
  F_l$.  We have the following:

  \begin{lem}\label{lem:data-C} 
    For $\chi\in \langle g_l\rangle^{\perp}$ denote by $M_{l,\chi}\inv
    $ the eigensheaf of $\OO_{\wt {E_l}}$ corresponding to $\chi$.
    Then the maps $\wt{E_l}\to E_{l, {a_l}}$ and $\wt{E_l}\to E_{l,
      b_l}$ induce isomorphisms:
$$L_{a_l, \chi}\inv \otimes\OO_{F_l}\cong M_{l,\chi}\inv (-A_{l,
  \chi}),\qquad L_{b_l, \chi}\inv\otimes\OO_{F_l}\cong M_{l, \chi}\inv
(-B_{l, \chi})$$ 
\end{lem}
\begin{proof} Follows by the normalization algorithm of \cite[\S
  3]{Pardini_AbelianCovers}.
\end{proof}
  Let $N_{l,\chi}:=A_{l,\chi}\cap B_{l,\chi}$ and let $T_{\chi}$ be the set  of points
  $y$ such that $C$ is singular at $y$
  and $\chi|_{H_y}$ is trivial.  
 We are now ready to describe  $(\im\alpha)_{ \chi}$:
   \begin{prop}\label{prop:Fchi} For every $\chi\in G^*\setminus \{1\}$,
  there is an exact sequence:
$$0\to (\im \alpha)_{\chi}\lra \oplus_{\{l|\chi(g_l)=1\}}M_{l,\chi}\inv(-N_{l,\chi}) \lra \OO_{T_{\chi}}\to 0.
$$
\end{prop}

\begin{proof}
In our setup, the map   $\wt{ B'}\to\wt {C'}$ is the disjoint union of two copies of $\wt B=\bigsqcup_{l=1}^s\wt{E_l}\to \bigsqcup_{l=1}^sF_l$ that are switched by the involution $j$. 
So by Lemma \ref{lem:data-C} the first sequence in   \eqref{eq:glue-surf} can be rewritten as:
\begin{equation}
   0\to \cF_{\chi}\to\oplus_{i=1}^t L_{i, \chi}\inv \to \oplus_{\{l|\chi(g_l)=1\}}M_{l,\chi}\inv.
 \end{equation}
In addition, if $F_l$ is a component of $C$ contained in $Y_{a_l}$ and $Y_{b_l}$,
 then again  by Lemma \ref{lem:data-C} the image of the map
 $L_{a_l, \chi}\inv\oplus L_{b_l,\chi}\inv \to M_{l,\chi}\inv$ is
 equal to $M_{l, \chi}\inv(-N^l_{\chi})$, so we have
 an exact sequence:
 \begin{equation}
   0\to(\im \alpha)_{\chi} \to \oplus_{\{l|\chi(g_l)=1\}}M_{l, \chi}\inv (-N_{\chi}^l) \to \cC_{\chi}\to 0,
 \end{equation}
where  the cokernel $\cC_{\chi}$ is concentrated on the set $T_{\chi}$.
 Using the description of the singularities of $X$ at these points
 given in \S \ref{ssec:setup-Z2}, one checks that $\cC_{\chi}$ has
 length 1 at points $y$ such that $\chi|_{H_y}$ is trivial and it is 0
 elsewhere, so $\cC_{\chi}=\OO_{T_{\chi}}$.
 \end{proof}

   \begin{remark}\label{rem:inertia}
     Let $y\in C$ be a smooth point, let $F$ be the irreducible
     component of $C$ that contains $y$ and let $Y_1$, $Y_2$ be the
     two components of $Y$ that contain $F$. Let $H$ the subgroup of
     $G$ generated by the inertia subgroups of $F$ and of the
     components of $D$ that contain $y$. Of course one has $H\subseteq
     H_y$, but in the present setup equality actually holds. Indeed,
     if $\chi\in H^{\perp}$ is a non trivial character, then by
     Proposition \ref{prop:Fchi} the second sequence in
     \eqref{eq:glue-surf} can be written near $y$ as $0\to\cF_{\chi}
     \to \OO_{Y_1}\oplus \OO_{Y_2}\overset{\alpha_{\chi}}{\to}
     \OO_F\to 0$, where $\alpha_{\chi}$ is given by $(f_1,f_2)\mapsto
     (f_1-f_2)|_F$. By Lemma \ref{lem:inertia}, there exist $z_i\in
     \OO_{Y_i}$, $i=1,2$, that correspond to functions on $X'_i$ that
     do not vanish at any point of $\pi\inv(y)$. Up to multiplying,
     say, $z_1$ by a nowhere vanishing regular function on $Y_1$ we
     can arrange that $z_{\chi}:=z_1-z_2\in \cF_{\chi}$. So $z_{\chi}$
     corresponds to a function on $X$ that is nonzero near
     $\pi\inv(y)$ and on which $G$ acts via the character $\chi$. It
     follows that $G/H$ acts freely on $\pi\inv(y)$, i.e. that
     $H=H_y$.
   \end{remark}

 We say that  a  point $y\in C$  is {\em relevant} iff either it is singular for $C$ or there exists $l$, $\chi$ with $\chi(g_l)=1$  such that $y\in N^l_{\chi}$.
Observe that, in view of the assumptions of \ref{ssec:setup-Z2}, by Proposition
  \ref{prop:index-graph} 
  and by the description of singularities of \S
  \ref{subsec:sings-red-base} the set of relevant points can be
  described intrinsically as the set of points of $C$ over which $X$
  is Gorenstein but not d.c..

\begin{cor}\label{cor:chi} Let ${\rm Rel}$ be the set of relevant
  points and let $\wt{B}=\sqcup_{l=1}^s \wt{E_l}$ be the normalization of the double locus $B$ of $X$.   Then:

$$\chi(\OO_X)=\chi(\OO_{X'})-\chi(\OO_{\wt{B}})+\sum_{y\in {\rm Rel}}[G:H_y].$$
\end{cor}
\begin{proof}
  Follows immediately by Proposition \ref{prop:Fchi}, by \eqref{eq:glue-surf} and by the fact that for $\chi=1$ one has the exact sequence:
  $$0\to (\im\al)_1\to \oplus_{l=1}^s\OO_{F_l}\to \OO_T\to 0,$$
  where $T$ is the set of singular points of $C$.
\end{proof}
We close this section by computing the numerical invariants of two of
the degenerations of Burniat surfaces described in
\cite{AlexeevPardini_CB}.
  \begin{example}\label{ex1} Let $G=\Z_2^2$, let $g_1, g_2,g_3$ be the
    nonzero elements of $G$ and for $i=1,2,3$ let $\chi_i\in G^*$ be
    the nonzero character such that $\chi_i(g_i)=1$.  Let
    $Y_1=\pp^1\times \pp^1$, $Y_2=\pp^2$ and let $Y$ be the surface
    obtained by gluing $Y_1$ and $Y_2$ along a smooth rational curve
    $C$ which is of type $(1,1)$ on $Y_1$ and is a line on $Y_2$. Fix
    three distinct points $y_1,y_2,y_3\in C$.  For $i=1,2,3$, let
    $D_{1,j}\subset Y_1$ be the union of a fibre and a section through
    $y_{j-1}$ and let $D_{2,j}\subset Y_2$ be a pair of lines through
    $y_{j+1}$ (the index $j$ varies in $\Z_3$). In the picture below
    $Y_1$ is represented on the left and $Y_2$ on the right, the curve
    $C$ is shown in green, red lines correspond to $D_{i,1}$, black
    lines to $D_{i,2}$ and blue lines to $D_{i,3}$.

    \begin{center}
      \includegraphics{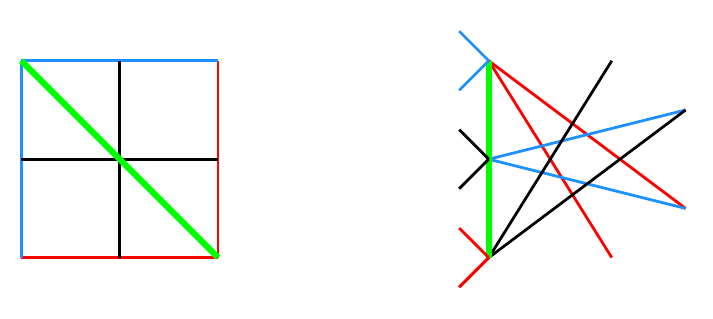}
    \end{center}

    For $i=1,2$, we let $\pi_i\colon X'_i\to Y_i$ be the standard
    $G$-cover with branch data $D_{i, j}$, $g_j$, $j=1,2,3$.
    Solving \eqref{eq:fundrel}, we get
    $L_{1,i}=\OO_{\pp^1\times\pp^1}(1,1)$ and $L_{2,j}=\OO_{\pp^2}(2)$,
    $j=1,2,3$, where $L_{i,j}\inv$ denotes the subsheaf of
    $\OO_{X'_i}$ corresponding to the character
    $\chi_j$. Notice that the line bundles $L_{i,j}\inv$ have no
    cohomology, hence in particular
    $\chi(\OO_{X'_1})=\chi(\OO_{X'_2})=1$.

    By \cite[\S 3]{Pardini_AbelianCovers}, for $i=1,2$ the
    normalization  of the cover of $C$ induced by $\pi_i$
    is the trivial  $G$-cover. So by Theorem
    \ref{thm:glue}, we can glue $X'_1\sqcup X'_2\to Y_1\sqcup Y_2$ to
    a cover $\pi\colon X\to Y$.  By \eqref{eq:canonical} we have:
$$K^2_X=4(K_{Y_1}+\frac{1}{2}(D_{1,1}+D_{1,2}+D_{1,3})+C)^2+4(K_{Y_2}+\frac{1}{2}(D_{2,1}+D_{2,2}+D_{2,3})+C)^2=2+4=6.$$
The curve $C$ is smooth and the points $y_1,y_2$ and $y_3$ are
relevant points with $H_{y_i}=G$, so Corollary \ref{cor:chi} gives:
$$\chi(\OO_X)=\chi(\OO_{X'_1})+\chi(\OO_{X'_2})-\chi(\OO_{\wt B})+[G:H_{y_1}]+[G:H_{y_2}]+[G:H_{y_3}]=1+1-4+1+1+1=1.$$
For $\chi=1$, we have an isomorphism 
$(\im\al)_1\cong  \OO_C$, hence $(\im\al)_1$ has no cohomology in degree $i>0$ and the exact sequence:
$$0\to\OO_Y\to \OO_{Y_1}\oplus\OO_{Y_2}\to (\im\al)_1=\OO_C\to 0$$
implies that $h^i(\OO_Y)=0$ for $i>0$. 
Next we compute the cohomology of the sheaves $\cF_{\chi}$.  By Proposition \ref{prop:Fchi}, for $j=1,2,3$ we have $(\im\al)_{\chi_j}=\OO_C(-y_j)$. So \eqref{eq:glue-surf} gives an exact 
 sequence:
$$0\to \cF_{\chi_j}\to L_{1,j}\inv \oplus L_{2,j}\inv \to \OO_C(-y_j)\to 0.$$
Therefore
$h^1(\cF_{\chi_j})=h^2(\cF_{\chi_j})=0$ for $j=1,2,3$ and thus $h^1(\OO_X)=h^2(\OO_X)=0$.
\end{example}

\begin{example}\label{ex:6cycle}
  Let $Y=Y_1\cup\dots \cup Y_6$ be the union of 6 copies of $\pp^2$
  glued in a cycle along lines as shown in the picture below.

  \begin{center}
  \includegraphics{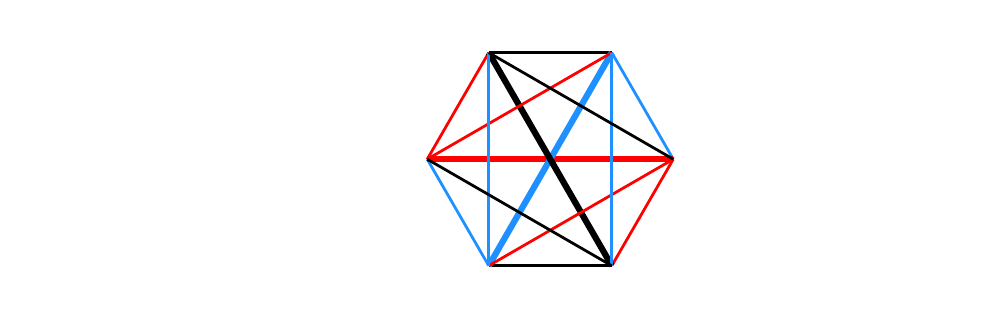}
  \end{center}
 
  As in the previous example, let $G=\Z_2^2$ and for $i\in\Z_6$ let
  $\pi_i\colon X'_i\to Y_i$ be the $G$-cover branched on the colored
  lines in the picture.  For every $i$, two of the sheaves
  $L_{i, \chi}$ are $\OO_{Y_1}(2)$ and the remaining one is
  $\OO_{Y_1}(1)$. So the $L_{i,\chi}\inv$ have no cohomology and
  $\chi(X'_i)=1$.  It's easy to check using Theorem \ref{thm:glue}
  that the cover $X'_1\sqcup \dots \sqcup X'_6\to Y_1\sqcup\dots
  \sqcup Y_6$ can be glued to a $G$-cover $\pi\colon X\to Y$.  The
  normalization $\wt B\to C$ of the induced cover of the double curve $C$ 
  is the disjoint union of 6 smooth rational curves, each mapping
  2-to-1 onto a component of $C$. The  only  relevant point is the singular point $y$ of
  $C$. So applying \eqref{eq:canonical} and Corollary
  \ref{cor:chi}, we get:
$$K^2_X=6,\quad \chi(\OO_X)=1.$$ Let $F_1,\dots F_6$ be the irreducible components of $C$. 
For $\chi=1$, as in the proof of Corollary \ref{cor:chi} we have an exact sequence:
$$0\to (\im\al)_1\to\oplus_{l=1}^6\OO_{F_l}\to \bK(y)\to 0,$$
which gives $h^i((\im\al)_1)=0$ for $i>0$. 
By Proposition \ref{prop:Fchi}, for $\chi\ne 0$ the sheaf $(\im\al)_{\chi}$ is isomorphic to the direct sum of two copies of $\OO_{\pp^1}$, hence it has no higher cohomology. So by \eqref{eq:glue-surf} we have $h^i(\cF_{\chi})=0$ for $i>0$ and therefore $h^1(\OO_X)=h^2(\OO_X)=0$.
\end{example}

\subsection{Singularities: the case $Y$
  smooth.}\label{subsec:sings-smooth-base}

We wish describe the singularities of a $\Z_2^r$-cover $\pi\colon X\to
Y$ as in \S \ref{ssec:setup-Z2}. Since the question is local, we fix
$y\in Y$ and we study $X$ locally above $Y$ in the \'etale topology.  By the assumptions in \S
\ref{ssec:setup-Z2}, the singularities of $X$ over a point $y\in Y$
lying on $q>2$ components of $Y$ are degenerate cusps such that the
exceptional divisor of its minimal semiresolution is a cycle of $q$
rational curves (cf.  \cite[def. 4.20]{KollarShepherdBarron}).  So it
is enough to analyze two cases:
\begin{itemize}
\item $Y$ is smooth,
\item $Y=Y_1\cup Y_2$ d.c. and $\pi$ is obtained by gluing standard
  covers $\pi_i\colon X'_i\to Y_i$, $i=1,2$.
\end{itemize}

\begin{remark} All the singularities listed in Tables
  ~\ref{tab:reduced-lines}--\ref{tab:R-D12non}, actually occur on some
  stable surface of general type.  To give examples of  the singularities  that appear  when the base  $Y$ of the cover is  smooth, one can
take
  $G=\Z_2^r$, $2\le r \le 4$, a set of generators $g_1,\dots g_k$ of
  $G$, $k\le 4$, and lines $L_1,\dots L_k$ through a point $y\in
  \pp^2$ such that the pair
$(\pp^2, (L_1+\dots L_k)/2)$
  is lc. 
  If $g=g_i$, define $D_{g_i}=L_i$, where $D'_i$ is a general curve of even degree and for $g\ne 1, g_1,\dots g_k$ let $D_g$ be a general curve of odd degree.
  The divisors  $D_g$ so defined  are the branch data for a
  $G$-cover $X\to\pp^2$ (equations \eqref{eq:fundrel} are easily seen
  have a solution in this case).   By Proposition \ref{prop:interesting-case} the surface $X$ is slc and it is of  general type as soon as the  the degree of the Hurwitz divisor $D$ is $>6$. There is only one point $x\in X$ that maps to $y$ and all the
  singularities listed in Tables 1, 2 and 3 with  can be
  realized as $(X,x)$ in this way and $|H|\ge 4$   (for the definition of $H$, see
  below).  The singularities with $|H|=2$ can be obtained by taking a
  double cover $X\to \pp^2$, branched on the sum of $k$ lines through
  $y$ and a general curve of degree $d$ such that $d+k$ is even and
  $\ge 8$. 
  
 Since  all the curves in the construction are general, the singularities of $X\setminus \{x\}$  are at most $A_1$ points.
  
  Similar constructions, slightly more involved, can be used to
  realize the singularities of
  Tables~\ref{tab:E-Dred}--\ref{tab:R-D12non}.
\end{remark}

We study the case $Y$ smooth in this section and the case $Y$
reducible in the next one.

We let $(D_1, g_1),\dots (D_k, g_k)$ be the branch data of $\pi$. We
may assume that $y\in D_i$ for every $i$. So by the condition that $D$
is slc we have $k\le 4$ and no three of the $D_i$ coincide. Whenever
the $D_i$ are not all distinct, we assume $D_1=D_2$.

All the possible cases are listed in
Tables~\ref{tab:reduced-lines},~\ref{tab:one-double-line},~\ref{tab:two-double-lines}
below.  The first digit in the label given to each case is equal to
the number $k$ of components through $y$, followed by $'$ if $D_1=D_2$
and by $''$ if $D_1=D_2$ and $D_3=D_4$ (obviously this case occurs
only for $k=4$).  
So, for instance, a label  of the form $3'.m$, where $m$ is any positive integer,  means that $y$ belongs to
three components of $D$, two of which coincide.

The entries in the columns have the following meaning:

\begin{itemize}
\item $|H|$: the order of the subgroup $H$ the subgroup generated by
  $g_1,\dots g_k$.
\item {\em Relations:} describes the relations between $g_1,\dots
  g_k$. For instance, $123$ means $g_1+g_2+g_3=0$.
\item{\em Singularity:} the notations are mostly
  standard. $\frac14(1,1)$ denotes a cyclic singularity $\bA^2/\bZ_4$
  with weights 1,1.  $T_{2,2,2,2}$ denotes an arrangement consisting
  of four disjoint $-2$-curves $G_1,\dots G_4$ and of a smooth
  rational curve $F$ intersecting each of the $G_i$ transversely at
  one point. The self intersection $F^2$ is given in the table.  In
  the non-normal case (Tables 2 and 3) we use the notations of
  \cite{KollarShepherdBarron}, where Koll\'ar and Shepherd-Barron
  classified all slc surface singularities over~$\bC$. We work in any
  characteristic $\ne2$ but only the singularities from the list in
  \cite{KollarShepherdBarron} appear.  ``deg.cusp$(k)$'' means a
  degenerate cusp (cf.  \cite[def. 4.20]{KollarShepherdBarron}) such
  that the exceptional divisor in the minimal semiresolution has $k$
  components.

\item $\iota$: the index of $x\in X$. It is equal to 1 if all the
  relations have even length and it is equal to 2 otherwise
  (cf. Proposition \ref{prop:index-normal}).

\item $\wX$: denotes the normalization of $X$ (the entries refer
  to the cases in Table~1).

\item $C_{\wX}\to C_X\to C_Y$: $C_{\wX}$ is the inverse image in
  ${\wX}$ of the double curve $C_X$ of $X$ and $C_Y$ is the image of
  $C_X$ in $Y$.  The symbol $\De$ denotes the germ of a smooth curve,
  and $\Ga_k$ is the seminormal curve obtained by gluing $k$ copies of
  $\De$ at one point.  The notation $\Ga_k\xrightarrow{a_1\dots a_k}C$
  means that the map restricts to a degree $a_i$ map on the $i$-th
  component of $\Ga_k$ (we do not specify the $a_i$ when they are all
  equal to 1).

\item $X\sr$: 
  is the minimal semiresolution of $X$. We write ``d.c.''
  when $X\sr$ has only normal crossings and ``pinch'' if it has also
  pinch points.
 \end{itemize}
\begin{theorem}
  The singularities of slc covers $\pi\colon X\to Y$ with smooth $Y$
  are listed in
  Tables~\ref{tab:reduced-lines},~\ref{tab:one-double-line},~\ref{tab:two-double-lines}.
\end{theorem}

\begin{table}[h]
  \centering
  \caption[t]{One, two, three, and four reduced lines}
  \begin{tabular}{|lllll|}
    \hline
    No. \Ttab &$|H|$ &Relations& $\iota$ &Singularity  \\

    \hline
    0.1 &1 &none& 1 &smooth \\
    \hline
    1.1 &2 &none& 1 &smooth\\

    \hline
    2.1 &4 &none& 1 &smooth \\
    2.2 &2 &12&1 &$A_1$ \\

    \hline
    3.1 &8 &none & 1& $A_1$\\

    3.2 &4 &12& 1 &$A_3$ \\

    3.3 &4 &123&2 &$\frac14(1,1)$ \\
    
    3.4 &2 &12,13&1 &$D_4$ \\

    \hline
    4.1 &16 &none & 1 &elliptic, $F^2=-4$\\
    4.2 &8 &12 & 1&elliptic, $F^2=-2$ \\
    4.3 &8 &123& 2 &$T_{2,2,2,2}$, $F^2=-4$\\
    4.4 &8 &1234 & 1&elliptic, $F^2=-8$ \\
    4.5 &4 &12 13 & 1&elliptic, $F^2=-1$ \\
    4.6 &4 &12 34 &1&elliptic, $F^2=-4$ \\
    4.7 &4 &12 134 & 2&$T_{2,2,2,2}$, $F^2=-3$ \\
    4.8 &2 &12 13 14 & 1&elliptic, $F^2=-2$ \\
    \hline
  \end{tabular}
  \label{tab:reduced-lines}
\end{table}

\begin{table}[h]
  \centering
  \caption[t]{Double line  + zero, one, or two reduced lines}
  \begin{tabular}{|llllllll|}
    \hline
    No.\Ttab &$|H|$ &Relations &$\iota$ &Singularity 
    &${\wX}$ &$C_{\wX}\to C_X\to C_Y$ &$X\sr$  \\

    \hline
    $2'.1$ &4 &none&1 &semismooth &2(1.1) &$2\De\to\De\to\De$ &d.c. \\
    $2'.2$ &2 &12 &1&semismooth  &$2(0.1)$ &$2\De\to\De\to\De$ &d.c.\\

    \hline
    $3'.1$ &8 &none&1 &semismooth  &2(2.1) &$2\De\to\De\xrightarrow{2}\De$&d.c.\\
    $3'.2$ &4 &12 &1&semismooth &2(1.1) &$2\De\to\De\xrightarrow{2}\De$ &d.c.\\
    $3'.3$ &4 &13 &1&semismooth  &(2.1)&$\De\xrightarrow{2}\De\to\De$  &pinch \\
    $3'.4$ &4 &123&2 &$(3'.1)/\bZ_2$ &2(2.2) &$2\De\to\De\to\De$ &d.c. \\
    $3'.5$ &2 &12 13 &1 &semismooth &(1.1) &$\De\xrightarrow{2}\De\to\De$  &pinch\\
    \hline

    $4'.1$ &16 &none&1 &deg.cusp$(2)$  &2(3.1) &$2\Ga_2 \to \Ga_2 \xrightarrow{22} \De$&d.c. \\
    $4'.2$ &8 &12 &1&deg.cusp$(2)$ &2(2.1) &$2\Ga_2 \to \Ga_2\xrightarrow{22} \De$ &d.c. \\
    $4'.3$ &8 &13 &1&deg.cusp$(1)$ &(3.1) &$\Ga_2\to \De\xrightarrow{2}\De$ &d.c. \\
    $4'.4$ &8 &34 &1 &deg.cusp(6)&2(3.2) & $2\Ga_2 \to \Ga_2\to\De$ &d.c. \\
    $4'.5$ &8 &123 &2 &$(4'.1)/\bZ_2$ &2(3.2) &$2\De\to \De\xrightarrow{2}\De$ &d.c. \\
    $4'.6$ &8 &134 &2 &$(4'.1)/\bZ_2$  &(3.1) &$\Ga_2\xrightarrow{22}\Ga_2\to\De$ &pinch\\
    $4'.7$ &8 &1234 &1&deg.cusp$(2)$ &2(3.3) &$2\Ga_2 \to \Ga_2 \to \De$ &d.c. \\
    $4'.8$ &4 &12 13&1 &deg.cusp$(1)$ &(2.1) &$\Ga_2\to\De\xrightarrow{2}\De$ &d.c.\\
    $4'.9$ &4 &13 14&1 &deg.cusp$(3)$ &(3.2) &$\Ga_2\to\De\to\De$&d.c.\\
    $4'.10$ &4 &12 34 &1&deg.cusp$(2)$ &2(2.2) &$2\Ga_2\to\Ga_2\xrightarrow{} \De$ &d.c.\\
    $4'.11$ &4 &13 24 &1&deg.cusp$(1)$ &(3.3) &$\Ga_2\xrightarrow{} \De\to\De$ &d.c.\\
    $4'.12$ &4 &12 134&2 &$(4'.2)/\bZ_2$&(2.1) &$\Ga_2\xrightarrow{22} \Ga_2\to\De$ &pinch \\
    $4'.13$ &4 &13 124 &2&$(4'.3)/\bZ_2$&(3.2) &$\De\xrightarrow{2}\De\to\De$ &pinch\\
    $4'.14$ &4 &123 34&2 &$(4'.4)/\bZ_2$&2(3.4) &$2\De\to \De\xrightarrow{}\De$ &d.c. \\
    $4'.15$ &2 &12 13 14&1 &deg.cusp$(1)$ & (2.2) &$\Ga_2 \to \De\to\De$  &d.c. \\

    \hline
  \end{tabular}
  \label{tab:one-double-line}
\end{table}

\vfill

\begin{table}[h]
  \centering
  \caption[t]{Two double lines}
  \begin{tabular}{|llllllll|}
    \hline
    No.\Ttab &$|H|$ &Relations &$\iota$ &Singularity 
    &${\wX}$ &$C_{\wX}\to C_X\to C_Y$ &$X\sr$ \\

    \hline
    $4''.1$\Ttab &16 &none& 1 &deg.cusp(4)& 4(2.1)& $4\Ga_2 \to \Ga_4\xrightarrow{2222}\Ga_2$&d.c. \\
    $4''.2$ &8 &12 &1&deg.cusp(4)&4(1.1)&$4\Ga_2\to \Ga_4\xrightarrow{2211} \Ga_2$&d.c. \\
    $4''.3$ &8 &13 &1 &deg.cusp(2)&2(2.1)&$2\Ga_2\to \Ga_2\xrightarrow{22}\Ga_2 $&d.c.\\
    $4''.4$ &8 &123&2 &$(4''.1)/\bZ_2$&2(2.1)&$2\Ga_2\xrightarrow{1122} \Ga_3\xrightarrow{211}\Ga_2$&pinch\\
    $4''.5$ &8 &1234 &1&deg.cusp(4)&4(2.2)&$4\Ga_2\to\Ga_4\to \Ga_2$&d.c.\\
    $4''.6$ &4 &12 13 &1&deg.cusp(2)&$2(1.1)$ &$2\Ga_2\to\Ga_2\xrightarrow{21}\Ga_2$& d.c.\\
    $4''.7$ &4 &12 34 &1&deg.cusp(4) &$4(0.1)$&$4\Ga_2\to \Ga_4\to\Ga_2$&d.c. \\
    $4''.8$ &4 &13 24 &1 &deg.cusp(2)&2(2.2)&$2\Ga_2\to\Ga_2\to\Ga_2$& d.c.\\
    $4''.9$ &4 &12 134 &2&$(4''.2)/\bZ_2$&2(1.1)&$2\Ga_2\xrightarrow{2211}\Ga_3\xrightarrow{}\Ga_2$&pinch \\
    $4''.10$ &4 &13 124 &2&$(4''.3)/\bZ_2$&(2.1)&$\Ga_2\xrightarrow{22}\Ga_2\to\Ga_2$&pinch \\
    $4''.11$ &2 &12 13 14 &1 &deg.cusp(2)&$2(0.1)$&$2\Ga_2\to\Ga_2\to\Ga_2$&d.c.\\

    \hline
  \end{tabular}
  \label{tab:two-double-lines}
\end{table}

Since all these singularities can be studied in a similar way, we just explain the method and work out two cases as an
illustration. We start with some general remarks:
\begin{itemize}
\item[(1)] we always assume $G=H$.  Indeed, the cover $\pi$ factors as
  $X\xrightarrow{\pi_2} X/H\xrightarrow{\pi_1} Y$.  By Lemma \ref{lem:inertia} the map $\pi_1$ is
  \'etale near $y$, while for every $z\in\pi_1\inv(y)$ the fiber
  $\pi_2\inv(z)$ consists only of one point.  Since $G$ acts
  transitively on each fiber of $\pi$, it is enough to describe the
  singularity of $X$ above any point $z\in\pi_1\inv(x)$.

\item[(2)] the cover $X$ is normal at $x$ iff $[D]=0$. It is nonsingular
  at $x$ iff either $k=1$ or $k=2$, $D_1\ne D_2$, $g_1\ne g_2$.
  Assume that $X$ is not normal, and let $F$ be an irreducible divisor
  that appears in $D$ with multiplicity 1. This means that, say,
  $F=D_1$ and $F=D_2$.  The normalization of $X$ along $F$ is a
  $G$-cover of $Y$ with branch data $(D_i,g_i)$, for $i\ne 1,2$, and,
  if $g_1+g_2\ne 0$, $(F, g_1+g_2)$ (cf. \cite[\S
  3]{Pardini_AbelianCovers}).

\item[(3)] the cover $X$ is said to be {\em simple} if the
  \underline{set} $\{g_1,\dots g_k\}$ is a basis of $|H|$ (for
  instance, $X$ is simple if the $g_i$ are all equal).  In this case
  $X$ is a complete intersection and it is very easy to write down 
  equations for it (see Case $4'.1$ below).

\item[(4)] the double curve $C_X$ maps onto the divisors that appear
  in $D$ with multiplicity $=1$. Since for a semismooth surface the
  double curve is locally irreducible, $X$ is never semismooth in the
  cases $4''$. In addition, if $X$ is semismooth then the pull back
  $C_{\wX}$ of $C_X$ to the normalization is smooth. Using this
  remark, it is easy to check that $X$ is never semismooth in the
  cases $4'$, either.

\item[(5)] in order to compute the minimal semiresolution $X\sr$, 
  we consider the blow up $\wh Y\to Y$ of $Y$ at $y$, pull back $X$ and
  normalize along the exceptional curve $E$ to get a cover $\wh X\to
  \wh Y$. The branch locus of $\wh X\to \wh Y$ is supported on a d.c. divisor and,
  by construction, the singularities of $\wh X$ are only of type $1$, $2$
  or $3'$. Looking at the tables, one sees that either $\wh X$ is
  semismooth or it has points of type $2.2$ or $3'.4$
  (cf. Table~\ref{tab:reduced-lines}). In the former case $\wh X$ is the
  minimal semiresolution. In the latter case, blowing up $\wh Y$ at the
  non semismooth points and taking base change and normalization along
  the exceptional divisor, one gets a semismooth cover $\wh{\wh X}\to
  \wh{\wh Y}$. The semiresolution $\wh{\wh X}\to X$ is minimal, except in cases
  $4''.5$, $4''.10$. In these cases the minimal semiresolution $X\sr$
  is obtained by contracting the inverse image in $\wh{\wh X}$ of the
  exceptional curve of the blow up $\wh Y\to Y$.
\end{itemize}

Next we analyze in detail two cases:

\underline{Case $4'.1$:} By remark (2) above, the normalization
${\wX}$ is an $H$-cover with branch data $(D_1, {g_1+g_2})$,
$(D_3,{g_3})$ and $(D_4, {g_4})$. So $g_1$ acts on $X$ without fixed
points and $X$ is the disjoint union of two copies of the cover
$(3.1)$.  We choose local parameters $u,v$ on $Y$ such that $D_1=D_2$
is given by $u=0$, $D_3$ is defined by $v=0$ and $D_4$ by $u+v=0$.
  
The cover $X$ is defined \'etale  locally above $y$ by the following equations:
\begin{equation}\label{eq4'.1}
  z_1^2=u, z_2^2=u, z_3^2=v, z_4^2=(u+v).
\end{equation}
 In particular $X$ is a complete
intersection (see remark (3) above). The element $g_i$ acts on $z_j$
as multiplication by $(-1)^{\de_{ij}}$. The double curve $C_X$ is the
inverse image of $u=0$, hence it is defined by $z_1=z_2=0, z_3=\pm
z_4$ and the map $C_X\to D_1$ is given by $z_3\mapsto z_3^2$, so $C_X$
is isomorphic to $\Ga_2$, with each component mapping $2$-to-$1$ to
$D_1\simeq \De$. The curve $C_{\wX}$ is the inverse image of $D_1$
in ${\wX}$, so it has two connected components, each isomorphic to
$\Ga_2$, that are glued together in the map ${\wX}\to X$.

To compute the minimal semiresolution, consider the blow up $\wh Y\to Y$
of $Y$ at $y$ and the cover $\wh X\to \wh Y$ obtained by pulling back $X\to
Y$ and normalizing along the exceptional curve $E$. The branch data
for $\wh X$ are $(E, g_1+g_2+g_3+g_4)$ and, for $i=1,\dots 4$, $(\wh{D_i},
g_i)$, where $\wh{}$ indicates the strict transform. The cover is singular
precisely above $\wh{ D_1}=\wh{ D_2}$ and it is easy, using the local equations,
to check that it is d.c. there. So $\wh X$ is the minimal semiresolution
of $X$. The exceptional divisor is the inverse image $F$ of $E$ in
$X$. Applying the normalization algorithm to the restricted cover
$F\to E$, one sees that the normalization $\wt F$ of $F$ is the
union of two smooth rational curves $F_1$ and $F_2$. The map
$\wt F\to F$ identifies the two points of $F_1$ that lie over the
point $E\cap D'_1$ with the corresponding two points of $F_2$. Hence
$\wh X$ is the minimal semiresolution of $X$ and the singularity is a
degenerate cusp solved by a cycle of two rational curves.  \smallskip

\underline{Case $4'.5$:} As in the previous case, ${\wX}$ and
$C_{\wX}$ can be computed by the normalization algorithm. One
obtains that ${\wX}$ is the disjoint union of two copies of $(3.2)$
and $C_{\wX}$ is the disjoint union of two copies of $\De$.  This
singularity is the quotient of a cover $X_0$ of type $(4'.1)$ by the
element $g_0:=g_1+g_2+g_3$. Since this element has odd length, the index
$\iota$ of $X$ at $x$ is equal to 2.

Since the only fixed point of $g_0$ on $X$ is $x:=\pi\inv (y)$, the
double curve $C_X$ is the quotient of the double curve $C_{X_0}$ of
$X_0$. The two components of $C_{X_0}$ are identified by $g_0$, thus
$C_X$ is irreducible and maps $2$-to-$1$ onto $D_1$.

To compute the minimal semiresolution, again we blow up $\wh Y\to Y$ at
$y$ and consider the cover $\wh X\to \wh Y$ obtained by pull back and
normalization along the exceptional curve $E$. As usual, we denote by
$\wh F$ the strict transform on $\wh Y$ of a curve $F$ of $Y$. The branch
data for $\wh X$ are $(\wh{D_1}, g_1)$, $(\wh{D_2}, g_2)$, $(\wh{D_3}, g_1+g_2)$,
$(\wh{D_4}, g_4)$, and $(E, g_4)$.  So $\wh X$ has normal crossings over
$\wh{D_1}$, it has four $A_1$ points over the point $\hat y:=\wh{D_4}\cap E$ and
it is smooth elsewhere (cf. Tables 1 and 2). We blow up at $\hat y$ and
take again pull back and normalization along the exceptional curve
$E_2$. We obtain a cover $\wh{\wh X}\to \wh{\wh Y}$ which is d.c. over the strict
transform $\wh{\wh{ D_1}}$ of $\wh {D_1}$ and has no other singularity, so $\wh{\wh X}\to
X$ is a semismooth resolution.  Let $E_1$ denote the strict transform
on $\wh{\wh Y}$ of the exceptional curve $E$ of the first blow up. Arguing as
in Case $4'.1$, one sees that inverse image of $E_1$ is the union of
two smooth rational curves $F_1^1$ and $F^1_2$ that intersect
transversely precisely at one point of the double curve, and the
inverse image of $E_2$ consists of 4 disjoint curves $F_2^1,\dots
F_2^4$. All these curves pull back to -2 curves on the normalization
of $\wh{\wh X}$ and, up to relabeling, $F_1^1, F^2_1, F^2_2$ and $F_1^2,
F^3_2, F^4_2$ form two disjoint $A_3$ configurations.  So $\wh{\wh {X}}$ is the
minimal semiresolution of $X$. In the notation of
\cite{KollarShepherdBarron}[def. 4.26], $\wh{\wh X}$ is obtained by gluing
two copies of $(A,\De)$ along $\De$.

\subsection{Singularities: the case $Y$
  reducible}\label{subsec:sings-red-base}

Here we repeat the local analysis of the previous section for the case
in which $Y=Y_1\cup Y_2$ is d.c., keeping as far as possible the same
notations.  So we fix $y\in C$, where $C$ is the double curve of $Y$,
and describe $X$ locally over $y$. We assume that $X\to Y$ is obtained
by gluing standard covers $\pi_i\colon X'_i\to Y_i$, $i=1,2$, such
that $y$ lies on all the components of the Hurwitz divisor $D$.  We
let $(D_1,g_1), \dots (D_k,g_k)$ be the union of the branch data of
$\pi_1$ and $\pi_2$ such that $D_i$ is distinct from the double curve
$C$ of $Y$ (hence $D=(D_1+\dots +D_k)/2$). We denote by $g_0$ the
generator of the inertia subgroup of $C$ for $\pi_1$ and $\pi_2$. By
Remark \ref{rem:inertia} the inertia subgroup $H_y$ is equal to $H:=
\langle g_0, g_1,\dots g_k\rangle$, so up to an \'etale cover we may
assume that $G=H$ and that $\pi\inv(y)=\{x\}$.

Since $D$ is $\Q$-Cartier, there are the same
number of $D_i$ on $Y_1$ and on $Y_2$. We order them so that all
components on $Y_1$ come first. Recall that $k\le 4$ by the assumption
that $(Y,D)$ is slc.    The cases in the tables are
labeled $E$ (``\'etale'') if $g_0=0$ and $R$ (``ramified'') if $g_0\ne
0$.  The first digit of the label is the number $k$ of branch lines
through $y$. It is followed by $'$ if $D_1=D_2$ and by $''$ if
$D_1=D_2$ and $D_3=D_4$.  For instance, in the cases E$4'.m$ the map $\pi$
is generically \'etale over $C$ and there are four branch lines
$D_1,\dots D_4$ with $D_1=D_2$, and $D_3\ne D_4$.

The singularities that we get here are non-normal, and as in
\cite[Thm.  4.21, 4.23]{KollarShepherdBarron} they turn out to be
either semismooth or degenerate cusps in the Gorenstein case and
$\Z_2$-quotients of these otherwise.

The tables here contain the same columns as those of \S
\ref{subsec:sings-smooth-base} plus an extra one, denoted $\chi$: this
is the contribution of $y$ in the formula for $\chi(\OO_X)$ of
Corollary \ref{cor:chi} (recall $|G|=2^r$). By Propositions  \ref{prop:index-non-normal}  and \ref{prop:index-graph}  the index $\iota$ is equal to 1 if all
relations have even length when reduced modulo $g_0$ and it is equal
to 2 otherwise.
\begin{theorem}
  The singularities of slc covers $\pi\colon X\to Y$ where $Y$ is the
  d.c. union of two smooth surfaces are given in
  Tables~\ref{tab:E-Dred}--\ref{tab:R-D12non}.
\end{theorem}

\begin{table}[h]
  \centering
  \caption[t]{$C$ not in the branch locus, zero, or two,  or four reduced lines}
\noindent\makebox[\textwidth]{%
  \begin{tabular}{|lllllllll|}
    \hline
    No.\Ttab &$|H|$ &Relations&$\iota$ &$\chi$&Singularity 
    &${\wX}$ &$C_{\wX}\to C_X\to C_Y$ &$X\sr$  \\

    \hline
    E$0.1$\Ttab &1&none&1&0 &d.c.&$(0.1)\sqcup (0.1)$ & $2\De\to \De\to\De $&d.c. \\
    \hline
    E$2.1$ &2&12 &1&0&d.c.&$(1.1)\sqcup (1.1)$ & $2\De\to \De\xrightarrow{2}\De $&d.c. \\
    \hline
    E$4.1$ &8&1234&1&$2^{r-3}$&deg.cusp$(4)$&$2(2.1)\sqcup 2(2.1)$&$2\Ga_2\sqcup2\Ga_2\to \Ga_4\xrightarrow{2222} \De$& d.c. \\
    E$4.2$ &4&12  34&1&$2^{r-2}$&deg.cusp$(4)$&$2(2.2)\sqcup 2(2.2)$ &$2\Ga_2\sqcup 2\Ga_2\to \Ga_4\to \De$& d.c. \\
    E$4.3$ &4&13  24&1&$2^{r-2}$&deg.cusp$(2)$&($2.1)\sqcup (2.1)$&$\Ga_2\sqcup\Ga_2\to \Ga_2\xrightarrow{22}\De$& d.c. \\
    E$4.4$ &2&12  13 14&1&$2^{r-1}$&deg.cusp$(2)$&$(2.2)\sqcup (2.2)$&$\Ga_2\sqcup\Ga_2\to \Ga_2\to\De$& d.c. \\
    \hline
  \end{tabular}
}
  \label{tab:E-Dred}
\end{table}

\begin{table}[h]
  \centering
  \caption[t]{$C$ not in the branch locus, a double line + two reduced lines.}
\noindent\makebox[\textwidth]{%
  \begin{tabular}{|lllllllll|}
    \hline
    No.\Ttab &$|H|$ &Relations&$\iota$& $\chi$ &Singularity 
    &${\wX}$ &$C_{\wX}\to C_X\to C_Y$ &$X\sr$ \\

    \hline
    E$4'.1$ \Ttab&8&1234 &1&$2^{r-3}$&deg.cusp$(6)$& $4(1.1)\sqcup 2(2.1)$&$4\Ga_2\sqcup 2\Ga_2\to \Ga_6\xrightarrow{112\dots 2} \Ga_2$&d.c. \\
    E$4'.2$ &4&12 34&1&$2^{r-2}$ &deg.cusp$(6)$&$4(0.1)\sqcup 2(2.2)$& $4\Ga_2\sqcup 2\Ga_2\to \Ga_6\to\Ga_2$& d.c.  \\
    E$4'.3$ &4&13 24 &1&$2^{r-2}$&deg.cusp$(3)$&$2(1.1)\sqcup (2.1)$&$2\Ga_2\sqcup \Ga_2\to \Ga_3\xrightarrow{122}\Ga_2$& d.c.  \\
    E$4'.4$ &2&12  13 14&1&$2^{r-1}$&deg.cusp$(3)$&$2(0.1)\sqcup (2.2)$ &$2\Ga_2\sqcup\Ga_2\to \Ga_3\to \Ga_2$& d.c. \\
    \hline
  \end{tabular}
}
  \label{tab:E-D1non-D2red}
\end{table}

\begin{table}[h]
  \centering
  \caption[t]{$C$ not in the branch locus, two pairs of double lines}
\noindent\makebox[\textwidth]{%
  \begin{tabular}{|lllllllll|}
    \hline
    No.\Ttab &$|H|$ &Relations  &$\iota$&$\chi$&Singularity
    &${\wX}$ &$C_{\wX}\to C_X\to C_Y$ &$X\sr$\\

    \hline
    E$4''.1$ \Ttab&8&1234&1 &$2^{r-3}$&deg.cusp$(8)$&$4(1.1)\sqcup 4(1.1)$& $4\Ga_2\sqcup 4\Ga_2\to \Ga_8\xrightarrow{112\dots211} \Ga_3$& d.c.\\
    E$4''.2$&4&12 34&1 &$2^{r-2}$&deg.cusp$(8)$&$4(0.1)\sqcup 4(0.1)$& $4\Ga_2\sqcup 4\Ga_2\to \Ga_8\to \Ga_3$& d.c. \\
    E$4''.3$ &4&13 24&1&$2^{r-2}$ &deg.cusp$(4)$&$2(1.1)\sqcup
    2(1.1)$&$2\Ga_2\sqcup 2\Ga_2\to \Ga_4\xrightarrow{1221}
\Ga_3$& d.c.\\
    E$4''.4$ &2&12  13 14&1&$2^{r-1}$&deg.cusp$(4)$&$2(0.1)\sqcup 2(0.1)$&$2\Ga_2\sqcup 2\Ga_2\to\Ga_4\to\Ga_3$& d.c. \\
    \hline

  \end{tabular}
}
  \label{tab:E-D12non}
\end{table}

\begin{table}[h]
  \centering
  \caption[t]{$C$ in the branch locus, zero, or two,  or four reduced lines}
\noindent\makebox[\textwidth]{%
  \begin{tabular}{|lllllllll|}
    \hline
    No.\Ttab &$|H|$ &Relations &$\iota$ &$\chi$&Singularity
    &${\wX}$ &$C_{\wX}\to C_X\to C_Y$ &$X\sr$  \\

    \hline
    R$0.1$\Ttab &2&none&1&0 &d.c.&$(1.1)\sqcup (1.1)$& $\De\sqcup \De\to\De\to\De$&d.c.\\
    \hline
    R$2.1$ &4&12&1&0&d.c.&$(2.1)\sqcup (2.1)$&$\De\sqcup\De\to\De\xrightarrow{2}{\De}$&d.c.\\
    R$2.3$ &2&12 01 &2&0&$({\rm R}2.1)/\Z_2$&$(2.2)\sqcup (2.2)$&$\De\sqcup\De\to\De\to\De$ &d.c.\\
    R$2.2$ &4&012&&& same as R$2.1$ &&&\\

    \hline
    R$4.1$ &16&1234&1&$2^{r-4}$ &deg.cusp$(4)$&$2(3.1)\sqcup 2(3.1)$&$2\Ga_2\sqcup 2\Ga_2\to\Ga_4\xrightarrow{2\dots 2}\De$&d.c.  \\
    R$4.2$ &8&1234 01&2&0&$({\rm R}4.1)/\Z_2$&$2(3.2)\sqcup(3.1)$&$2\De\sqcup \Ga_2\to \Ga_2\xrightarrow{22}\De$ & d.c.\\
    R$4.3$ &8&1234 012 &1&$2^{r-3}$ &deg.cusp$(4)$&$2(3.3)\sqcup 2(3.3)$&$2\Ga_2\sqcup2\Ga_2\to\Ga_4\to\De$& d.c.  \\ 
    R$4.4$ &8&1234 013&1&$2^{r-3}$&deg.cusp$(2)$&$(3.1)\sqcup (3.1)$&$\Ga_2\sqcup\Ga_2\to\Ga_2\xrightarrow{22}\De$& d.c. \\
    R$4.5$ &8&12  34&1&$2^{r-3}$&deg.cusp$(12)$&$2(3.2)\sqcup 2(3.2)$&$2\Ga_2\sqcup2\Ga_2\to\Ga_4\to\De$& d.c. \\
    R$4.6$ &4&12  34 01&2&0 &$(\rm{R}4.5)/\bZ_2$&$2(3.4)\sqcup (3.2)$& $2\De\sqcup \Ga_2\to \Ga_2\to \De$&d.c. \\
    R$4.7$ &4&12 34 013 &1&$2^{r-2}$&deg.cusp$(6)$&$(3.2)\sqcup(3.2)$& $\Ga_2\sqcup\Ga_2\to\Ga_2\to\De$ & d.c. \\
    R$4.8$&8&13 24&&& same as R$4.4$ & & &  \\
    R$4.9$ &4&13 24 01&2&0&$(\rm{R}4.8)/\bZ_2$&$(3.2)\sqcup (3.2)$&$\De\sqcup\De\to\De\xrightarrow{2}\De$& d.c. \\
    R$4.10$ &4&   13 24 012 &1&$2^{r-2}$&deg.cusp$(2)$&$(3.3)\sqcup(3.3)$&$\Ga_2\sqcup\Ga_2\to\Ga_2\to\De$& d.c. \\
    R$4.11$&4&12 13 14&&& same as R$4.7$ & & &   \\
    R$4.12$ &2&12  13 14 01&2&0&$(\rm{R}4.11)/\bZ_2$&$(3.4)\sqcup (3.4)$&$\De\sqcup \De\to\De$& d.c. \\
    R$4.13$&16&01234& &&same as R$4.1$ & & &   \\
    R$4.14$ &8&12 034&1&$2^{r-3}$&deg.cusp$(8)$&$2(3.2)\sqcup 2(3.3)$&$2\Ga_2\sqcup2\Ga_2\to\Ga_4\to\De$& d.c. \\
    R$4.15$&8&13 024&&& same as R$4.4$ & & &  \\
    R$4.16$&8&123 04&&& same as R$4.2$ & & &   \\
    R$4.17$&4&12 13 014&1&$2^{r-2}$&deg.cusp$(4)$&$(3.2)\sqcup(3.3)$&$\Ga_2\sqcup\Ga_2\to\Ga_2\to\De$& d.c. \\ 
    R$4.18$ &4&12  134 01&2&0&$(\rm{R}4.14)/\bZ_2$&$2(3.4)\sqcup (3.3)$&$2\De\sqcup \Ga_2\to \Ga_2\to \De$& d.c.\\
    R$4.19$&4&13 124 01& &&same as R$4.9$& & &  \\
    \hline

  \end{tabular}
}
  \label{tab:R-Dred}
\end{table}
\newpage

\begin{table}[h]
  \centering
  \caption[t]{$C$ in the branch locus, a double line + two reduced lines.}

\noindent\makebox[\textwidth]{%
  \begin{tabular}{|lllllllll|}
    \hline
    No.\Ttab &$|H|$ &Relations &$\iota$&$\chi$&Singularity 
    &${\wX}$ &$C_{\wX}\to C_X\to C_Y$ &$X\sr$  \\
    \hline
    R$4'.1$& 16 & 1234&1&$2^{r-4}$ & deg.cusp$(6)$ &$4(2.1)\sqcup 2(3.1)$ & $4\Ga_2\sqcup2\Ga_2\to\Ga_6 \xrightarrow{2\dots 2}\Ga_2$ &d.c.  \\
    R$4'.2$& 8 & 1234 01&2&0 &$(\rm{R}4'.1)/\Z_2$ &$2(2.1)\sqcup (3.1)$  &$2\Ga_2\sqcup\Ga_2\xrightarrow{221111}\Ga_4 \xrightarrow{1122}\Ga_2$  & pinch  \\
    R$4'.3$& 8 & 1234 03& 2&0& $(\rm{R}4'.1)/\Z_2$&$2(2.1)\sqcup 2(3.2)$ &$2\Ga_2\sqcup 2\De\to \Ga_3\xrightarrow{222}\Ga_2$&d.c. \\
    R$4'.4$& 8 & 1234 012&1&$2^{r-3}$&deg.cusp$(6)$&$4(2.2)\sqcup 2(3.3)$&$4\Ga_2\sqcup2\Ga_2\to \Ga_6\to\Ga_2$& d.c. \\ 
    R$4'.5$& 8 & 1234 013&1&$2^{r-3}$&deg.cusp$(3)$&$2(2.1)\sqcup (3.1)$&$2\Ga_2\sqcup \Ga_2\to \Ga_3\xrightarrow{222}\Ga_2$& d.c.\\

    R$4'.6$ &8&12  34&1&$2^{r-3}$&deg.cusp$(10)$&$4(1.1)\sqcup 2(3.2)$&$4\Ga_2\sqcup2\Ga_2\to \Ga_6\xrightarrow{221\dots 1}\Ga_2$& d.c.\\
    R$4'.7$& 4 & 12 34 01&2&0&$(\rm{R}4'.6)/\Z_2$&$2(1.1)\sqcup(3.2)$&$2\Ga_2\sqcup\Ga_2\xrightarrow{221\dots1} \Ga_4\to\Ga_2$& pinch \\

    R$4'.8$& 4 & 12 34 03&2&0&$(\rm R4'.6)/\Z_2$&$2(1.1)\sqcup2(3.4)$&$2\Ga_2\sqcup2\De\to\Ga_3\xrightarrow{211}\Ga_2$& d.c. \\
 
    R$4'.9$& 4 & 12 34 013&1&$2^{r-2}$&deg.cusp$(5)$&$2(1.1)\sqcup(3.2)$&$2\Ga_2\sqcup\Ga_2\to\Ga_3\xrightarrow{211}\Ga_2$&d.c. \\
    R$4'.10$& 8 & 13 24 & & &same as R$4'.5$ & & &\\
    R$4'.11$& 4 & 13 24 01&2&0 &$(\rm R4'.10)/\Z_2$ &$(2.1)\sqcup (3.2)$ &$\Ga_2\sqcup \De\xrightarrow{211}\Ga_2\xrightarrow{12} \Ga_2$ &\\
    R$4'.12$& 4 & 13 24 012&1&$2^{r-2}$&deg.cusp$(3)$&$2(2.2)\sqcup(3.3) $&$2\Ga_2\sqcup\Ga_2\to\Ga_3\to\Ga_2$&d.c.\\
    R$4'.13$& 4 & 12 13 14& & &same as R$4'.9$ & & &\\
    R$4'.14$& 2 & 12 13 14 01&2& 0&$(\rm R4'.13)/\Z_2$&$(1.1)\sqcup(3.4)$&$\Ga_2\sqcup\De\xrightarrow{211}\Ga_2\to\Ga_2$& pinch \\
    R$4'.15$& 8 &13 024 &&&same as R$4'.5$  & & &\\
    R$4'.16$& 8 & 12 034 &1&$2^{r-3}$ &deg.cusp$(6)$&$4(1.1)\sqcup 2(3.3)$&$4\Ga_2\sqcup2\Ga_2\to \Ga_6\xrightarrow{221\dots 1}\Ga_2$& d.c. \\
    R$4'.17$& 8 &13 024& &&same as R$4'.5$  & & &\\
    R$4'.18$& 8 & 34 012&1&$2^{r-3}$ &deg.cusp$(10)$&$4(2.2)\sqcup 2(3.2)$&$4\Ga_2\sqcup2\Ga_2\to \Ga_6\to\Ga_2$& d.c. \\
    R$4'.19$& 8 &123 04& && same as R$4'.3$  & & &\\
    R$4'.20$& 8 &134 02& & & same as R$4'.2$  & & &\\
    R$4'.21$& 4 & 12 13 014&1&$2^{r-2}$ &deg.cusp$(3)$&$2(1.1)\sqcup(3.3)$&$2\Ga_2\sqcup\Ga_2\to\Ga_3\xrightarrow{211}\Ga_2$& d.c. \\
    R$4'.22$& 4 & 13 14 012&1&$2^{r-2}$ &deg.cusp$(5)$&$2(2.2)\sqcup(3.2)$&$2\Ga_2\sqcup\Ga_2\to\Ga_3\to\Ga_2$& d.c.\\
    R$4'.23$& 4 & 12 134 01 &2&0&$(\rm{R}4'.16)/\Z_2$&$2(1.1)\sqcup(3.3)$&$2\Ga_2\sqcup \Ga_2\to \Ga_3\xrightarrow{211}\Ga_2$& pinch \\

    R$4'.24$& 4 & 13 124 01& & & same as R$4'.11$ & & &\\
    R$4'.25$& 4 & 34 123 03 &2&0& $(\rm R4'.18)/\Z_2$&$2(2.2)\sqcup2(3.4)$&$2\Ga_2\sqcup 2\De\to\Ga_3\to\Ga_2$& d.c. \\
    \hline
  \end{tabular}
}
  \label{tab:R-D1non-D2red}
\end{table}
\newpage

\begin{table}[h]
   \centering
  \caption[t]{$C$ in the branch locus, two pairs of double lines}    
\noindent\makebox[\textwidth]{%
  \begin{tabular}{|lllllllll|}
    \hline
    No.\Ttab &$|H|$ &Relations&$\iota$&$\chi$ &Singularity 
    &${\wX}$ &$C_{\wX}\to C_X\to C_Y$ &$X\sr$  \\
    \hline
    R$4''.1$ &16&1234&1&$2^{r-4}$& deg.cusp$(8)$&$4(2.1)\sqcup 4(2.1)$&$4\Ga_2\sqcup4\Ga_2\to\Ga_8\xrightarrow{2\dots 2}\Ga_3$&d.c.  \\
    R$4''.2$ &8&1234 01&2&0&$({\rm R}4''.1)/\Z_2$&$2(2.1)\sqcup 2(2.1)$&$2\Ga_2\sqcup2\Ga_2\xrightarrow{221\dots 1}\Ga_5\xrightarrow{11222}\Ga_3$& pinch \\
    R$4''.3$ &8&1234 012&1&$2^{r-3}$ &deg.cusp$(8)$&$4(2.2)\sqcup4(2.2)$&$4\Ga_2\sqcup4\Ga_2\to\Ga_8\to\Ga_3$& d.c. \\
    R$4''.4$ &8&1234 013&1&$2^{r-3}$&deg.cusp$(4)$&$2(2.1)\sqcup 2(2.1)$&$2\Ga_2\sqcup2\Ga_2\to \Ga_4\xrightarrow{2222}\Ga_3$& d.c. \\
    R$4''.5$ &8&12  34&1&$2^{r-3}$&deg.cusp$(8)$&$4(1.1)\sqcup 4(1.1)$&$4\Ga_2\sqcup 4\Ga_2\to \Ga_8\xrightarrow{22111122}\Ga_3$& d.c. \\
    R$4''.6$ &4&12  34 01 &2&0 &$(\rm R 4''.5)/\Z_2$&$2(1.1)\sqcup 2(1.1)$&$2\Ga_2\sqcup2\Ga_2\xrightarrow{221\dots1}\Ga_5\xrightarrow{11112} \Ga_3$& pinch \\
    R$4''.7$ &4&12 34 013 &1&$2^{r-2}$&deg.cusp$(4)$&$2(1.1)\sqcup 2(1.1)$&$2\Ga_2\sqcup2\Ga_2\to\Ga_4\xrightarrow{2112}\Ga_3$& d.c. \\
    R$4''.8$&8&13 24& & & same as R$4''.4$ & & &  \\
    R$4''.9$ &4&13 24 01&2&0 &$({\rm R}4''.8)/\Z_2$&$(2.1)\sqcup(2.1)$&$\Ga_2\sqcup\Ga_2\xrightarrow{212} \Ga_3\xrightarrow{121}\Ga_3$& pinch\\
    R$4''.10$ &4&   13 24 012&1&$2^{r-2}$ &deg.cusp$(4)$&$2(2.2)\sqcup 2(2.2)$&$2\Ga_2\sqcup2\Ga_2\to\Ga_4\to\Ga_3$& d.c. \\

    R$4''.11$&4&12 13 14& &&same as R$4''.7$&  & &  \\
    R$4''.12$ &2&12  13 14 01&2&0&$(\rm R4''.11)/\Z_2$&$(1.1)\sqcup (1.1)$&$\Ga_2\sqcup\Ga_2\xrightarrow{2112} \Ga_3\to\Ga_3$& pinch \\

    R$4''.13$&16&01234& & & same as R$4''.1$ & & &  \\
    R$4''.14$ &8&12 034&1& $2^{r-3}$ &deg.cusp$(8)$&$4(1.1)\sqcup4(2.2)$&$4\Ga_2\sqcup4\Ga_2\to\Ga_8\xrightarrow{221\dots1}\Ga_3$& d.c. \\

    R$4''.15$&8&13 024& && same as R$4''.4$ & & &   \\
    R$4''.16$&8&123 04& && same as R$4''.2$ & & &   \\
    R$4''.17$&4&12 13 014&1& $2^{r-2}$ &deg.cusp$(4)$&$2(1.1)\sqcup 2(2.2)$&$2\Ga_2\sqcup 2\Ga_2\to\Ga_4\xrightarrow{2111}\Ga_3$& d.c. \\
    R$4''.18$ &4&12  134 01&2&0&$(\rm R 4''.14)/\Z_2$&$2(1.1)\sqcup 2(2.2)$&$2\Ga_2\sqcup 2\Ga_2\xrightarrow{221\dots1} \Ga_5\to\Ga_3$& pinch \\

    R$4''.19$&4&13 124 01& &&same as R$4''.9$ & & &  \\
    \hline

  \end{tabular}
}
  \label{tab:R-D12non}
\end{table}
  
The analysis of the singularities in the reducible case is similar to
the case $Y$ smooth. One blows up $Y$ at the point $y$ and takes pull
back and normalization of $X$ along the exceptional divisor. Repeating
this process, if necessary, one obtains a semiresolution $X_0\to X$.
If $X_0$ is not minimal, then the minimal semiresolution $X\sr\to X$ is
obtained by blowing down the $-1$-curves of $X_0$.

As the computations are all similar,  we work out a only a couple of cases to
show the method.  \smallskip

\underline{Case R$4'.1$:} the normalization $\wt X$ is equal to
$\wt{ X'_1}\sqcup \wt{X'_2}$, where $\wt{ X'_i}$ is the normalization of
$X'_i$. The branch data of $\wt{ X'_1}\to Y_1$ are $(D_1, g_1+g_2)$,
$(D_0,g_0)$, so $\wt{ X'_1}$ is \'etale locally  the disjoint union of four
copies of the cover $(2.1)$. $X'_2=\wt{X'_2}$ is  \'etale locally  the disjoint union of two copies  
$(3.1)$.

The image $C_Y$ of the double curve $C_X$ is equal to $C\cup D_1$. The
preimage in $\wt{ X'_1}$ of $C_Y$ is the disjoint union of four copies of
$\Ga_2$. The preimage of $C_Y$ in $\wt{X'_2}$ is equal to two copies of
$\Ga_2$. Hence $C_{\wt X}=4\Ga_2\sqcup 2\Ga_2$. Each component of
$C_{\wt X}$ maps 2-to-1 onto its image. The map $C_{\wt X}\to C_X$
identifies in pairs the four components of the preimage of $D_1$ and
the eight components of the preimage of $C$. So $C_X$ is $\Ga_6$, with
two components mapping 2-to-1 onto $D_1$ and four components mapping
2-to-1 onto $C$.

To compute the semiresolution, blow up $y\in Y$ to get $\wh Y\to Y$. Let
$E_1\subset Y_1$ and $E_2\subset Y_2$ be the irreducible components of
the exceptional divisor. Let $\wh \pi\colon \wh X\to \wh Y$ be the $G$-cover
obtained from $X\to Y$ by taking pull back and normalizing along $E_1$
and $E_2$. Denoting by $\wh{}$ the strict transform on $\wh Y$, the branch
data of $\wh \pi$ are: $(E_1, g_0+g_1+g_2)$, $(E_2,
g_0+g_3+g_4=g_0+g_1+g_2)$, $(\wh{D_1}, g_1)$, $(\wh{D_2}=\wh{D_1}, g_2)$,
$(\wh{D_3},g_3)$, $(\wh{D_4},g_4)$, $(\wh C, g_0)$.  So $\wh X$ is d.c. by the
tables of \S \ref{subsec:sings-smooth-base} and it is therefore the
semiresolution $X\sr$ of $X$. The preimage of $E_1$ is the union of
four smooth rational curves meeting in pairs over the point $E_1\cap
\wh{D_1}$. The preimage of $E_2$ is the disjoint union of two rational
curves, which together with the components of the preimage of $E_1$
form a cycle of six rational curves. The singularity $x\in X$ is
Gorenstein by Proposition  \ref{prop:index-graph}   
hence it is ``deg.cusp$(6)$''.  \smallskip

\underline{Case R$4'.2$:} This is a $\Z_2$-quotient of R$4'.2$ and it
is not Gorenstein by Proposition \ref{prop:index-non-normal}. The
normalization $\wt X$ is equal to $\wt{ X'_1}\sqcup \wt{X'_2}$, where
$ \wt{X'_i}$ is the normalization of $ X'_i$. The branch data of
$\wt{ X'_1}\to Y$ are $(D_1, g_0+g_2)$, $(D_0,g_0)$, so $\wt{ X'_1}$ is \'etale locally
the disjoint union of two copies of the cover $(2.1)$.  The image
$C_Y$ of the double curve $C_X$ is equal to $C\cup D_1$. The preimage
in $\wt{ X_1}$ of $C_Y$ is the disjoint union of two copies of $\Ga_2$.
The preimage of $C_Y$ in $\wt{X'_2}$ is $\Ga_2$. Hence
$C_{\wh X}=2\Ga_2\sqcup \Ga_2$.  Each component of $C_{\wt X}$ maps 2-to-1
onto its image in $C_Y$.  The map $C_{\wt X}\to C_X$ identifies glues two itself each of the  two
components of the preimage of $D_1$ and it identifies in pairs the
four components of the preimage of $C$.  So $C_X$ is $\Ga_4$, with two
components mapping 1-to-1 onto $D_1$ and two components mapping 2-to-1
onto $C$.

To compute the semiresolution, blow up $y\in Y$ to get $\wh Y\to Y$. Let
$E_1\subset Y_1$ and $E_2\subset Y_2$ be the irreducible components of
the exceptional divisor. Let $\wh \pi\colon \wh X\to \wh Y$ be the $G$-cover
obtained from $X\to Y$ by taking pull back and normalizing along $E_1$
and $E_2$.  Denoting by $\wh{}$ the strict transform on $\wh Y$, the branch
data of $\wh \pi$ are: $(E_1, g_2)$, $(E_2, g_0+g_3+g_4=g_2)$, $(\wh{D_1},
g_1=g_0)$, $(\wh{D_2}=\wh{D_1}, g_2)$, $(\wh{D_3},g_3)$, $(\wh{D_4},g_4)$, $(\wh C,
g_0)$. By the tables of \S \ref{subsec:sings-smooth-base}, $\wh X$ has
two pinch points over the point $\wh{D_1}\cap E_1$ and is at most
d.c. elsewhere, hence it is equal to the minimal semiresolution $
X\sr$.  The preimage of $E_1$ is a pair of smooth rational curves meeting
over the point $E_1\cap \wh{D_1}$. The preimage of $E_2$ is a smooth
rational curve, meeting each component of the preimage of $E_1$ at a
point lying over $\wh C\cap E_1=\wh C\cap E_2$.

In the notation of \cite{KollarShepherdBarron}[Def. 4.26], $X\sr$ is a
chain consisting of copy of $(A,2\De)$ (namely the second component of
$X\sr$) in the middle and two copies of $(A, 2\De)$ with $\De$ pinched
at the ends.

\bibliographystyle{amsalpha}
\renewcommand\MR[1]{}
\def\cprime{$'$}
\providecommand{\bysame}{\leavevmode\hbox to3em{\hrulefill}\thinspace}
\providecommand{\MR}{\relax\ifhmode\unskip\space\fi MR }
\providecommand{\MRhref}[2]{%
  \href{http://www.ams.org/mathscinet-getitem?mr=#1}{#2}
}
\providecommand{\href}[2]{#2}

\end{document}